\documentclass[%
  onecolumn
   , hidempi
]{mpi2015-cscpreprint}

\usepackage[american]{babel}
\usepackage{subfig}

\usepackage{graphicx}
\usepackage{pgfplots}
\usepackage{tikz}

\usepackage[textwidth=2.5cm, textsize=small]{todonotes}

\usepackage{amssymb}
\usepackage{amsthm}
\usepackage{algorithm}
\usepackage{algorithmic}
\renewcommand{\algorithmicrequire}{\textbf{Input:}}

\usepackage{booktabs,dcolumn,array}
\newcolumntype{d}{D{.}{.}{2.3}}
\newcolumntype{C}{>{\centering}p}

\newtheorem{theorem}{Theorem}[section]

\newtheorem{remark}{Remark}

\crefname{appsec}{Appendix}{Appendices}
\usepackage{appendix}

\newlength{\fheight}
\newlength{\fwidth}

\input{mymymacros.sty}

\begin{document}


\title{Accurate error estimation for model reduction of nonlinear dynamical systems via data-enhanced error closure}

\author[$\ast$]{Sridhar Chellappa}
\affil[$\ast$]{Max Planck Institute for Dynamics of Complex Technical Systems, 39106 Magdeburg, Germany.\authorcr
\vspace{1mm}
	\email{chellappa@mpi-magdeburg.mpg.de}, \orcid{0000-0002-7288-3880} \authorcr
	\email{feng@mpi-magdeburg.mpg.de}, \orcid{0000-0002-1885-3269} \authorcr
	\email{benner@mpi-magdeburg.mpg.de}, \orcid{0000-0003-3362-4103}}
  
\author[$\ast$]{Lihong Feng}

\author[$\ast$]{Peter Benner}
  
\shorttitle{Data-enhanced error closure for MOR}
\shortauthor{S. Chellappa, L. Feng, P. Benner}
\shortdate{}
  

  
\abstract{%
Accurate error estimation is crucial in model order reduction, both to obtain small reduced-order models and to certify their accuracy when deployed in downstream applications such as digital twins. In existing a posteriori error estimation approaches, knowledge about the time integration scheme is mandatory, e.g., the residual-based error estimators proposed for the reduced basis method. This poses a challenge when automatic ordinary differential equation solver libraries are used to perform the time integration. To address this, we present a data-enhanced approach for \emph{a posteriori} error estimation. Our new formulation enables residual-based error estimators to be independent of any time integration method. To achieve this, we introduce a corrected reduced-order model which takes into account a data-driven closure term for improved accuracy. The closure term, subject to mild assumptions, is related to the local truncation error of the corresponding time integration scheme. We propose efficient computational schemes for approximating the closure term, at the cost of a modest amount of training data. Furthermore, the new error estimator is incorporated within a greedy process to obtain parametric reduced-order models. Numerical results on three different systems show the accuracy of the proposed error estimation approach and its ability to produce ROMs that generalize well.%
}

 \novelty{\begin{itemize}
 		\item We present a new \apost~error estimation framework for nonlinear dynamical systems, which is independent of the ordinary differential equation solver used.
 		
 		\item Our approach imposes a simpler, user-defined time integration scheme on the high-fidelity snapshots data and learns a parametrized closure term to recover the true residual needed for error estimation.
 		
 		\item We numerically demonstrate the approximability of the data-driven closure term and propose efficient algorithms for its computation.
 		
 		\item The proposed methodology eliminates one of the main intrusive elements present in the reduced basis method, viz., the knowledge of the time integration scheme; it further paves way for the reduced basis method to become better integrated within existing computational software packages such as \MATLAB $\,$ or Python.
 		
 \end{itemize}}

\maketitle

  
\section{Introduction}%
\label{sec:intro}

Model order reduction (MOR) has become an important enabling technology that facilitates the rapid and reliable simulations of large-scale systems in a number of scientific disciplines. The central goal of MOR is to replace an \emph{expensive-to-compute} full-order model (FOM) with a surrogate, called the reduced-order model (ROM). The ROM has fewer degrees of freedom which enables real-time simulation -- an important requirement in recent developments such as \emph{digital twins}~\cite{morharthw18,morKapetal22}. We refer to the recent books \cite{morBenetala21,morBenetalb21,morBenetalc21} for a detailed background and the state-of-the-art in MOR literature. 

In this work, we are interested in obtaining ROMs for parametric, nonlinear dynamical systems that arise from the numerical discretization of partial differential equations (PDEs). For such systems, the Reduced Basis Method (RBM)~\cite{morRozHP08,morQuaMN16,morHesRS16} is a commonly used MOR approach. The RBM relies on an \apost~error estimator to perform an iterative, greedy sampling of the parameter space to compute solution snapshots, in order to build a linear subspace that serves as an approximation for the solution space. The FOM equations are then projected onto this subspace to obtain the final ROM. \emph{A posteriori}, residual-based error bounds/estimators are an indispensable part of the RBM; such estimators have been proposed for a variety of system classes such as coercive and non-coercive elliptic systems~\cite{morMacMP01,morRov03}, linear parabolic systems~\cite{morGreP05}, non-linear or non-affine systems~\cite{morVeretal03,morCanTU09,morGre12}. While initial development in the RBM community focused on deriving \apost~error estimators in the variational setting, more recent work has also focused on deriving such estimators for PDEs that are discretized using finite volume method or the finite difference method; see \cite{morHaaO11,morDroHO12,morWirSH14,morZhaFLetal15,morCheFB19a}. 

\subsection{Motivation}
For dynamical systems, all of the above works on \apost~error estimation assume prior knowledge of the time integration scheme. The expression of the time integration method is needed to compute the residual vector, which is obtained by plugging an approximate solution in the residual expression. For instance, if a Runge-Kutta method were used, the Butcher tableau of the corresponding scheme needs to be known. However, there are scenarios where no (or only incomplete) knowledge of the time integration scheme is available. One such example is when automatic ordinary differential equation (ODE) solver libraries, e.g., SUNDIALS~\cite{gardner2022sundials,hindmarsh2005sundials}, ODEPACK~\cite{hindmarsh1983odepack}, ARKODE~\cite{reynolds2022arkode} are used to perform the time integration. Such ODE solvers are often available readily in computational software such as \MATLAB $\,$ (\texttt{ode45, ode15s}, etc.) or Python (\texttt{scipy.integrate.odeint}). Another popular ODE solver library is the \textsf{TS} library~\cite{AbhyankarEtAl2018} available in PETSc~\cite{petsc-web-page}. Typically, automatic ODE solver packages are computationally efficient as they employ adaptive order selection (e.g., \texttt{ode15s} can switch between order 1 to order 5 methods) and/or adaptive choice of the time step, which is beneficial in case the system being considered exhibits \emph{stiffness}. Moreover, much time is saved by re-using efficient and robust software which is already available and well-maintained. Whenever automatic ODE solvers are used, the standard RBM can no longer be used in a straightforward manner to obtain a ROM. This is owing to the fact that the exact expression of the time integration method used within such ODE solvers is unknown and deriving the corresponding error estimator is an open question. This serves as the fundamental motivation for our work.

\subsection{Main contributions}
Our work introduces a new data-enhanced paradigm for \apost, residual-based error estimation for the RBM, in the absence of any knowledge about the time integration scheme used. In addition to the discretized system matrices, our approach only requires access to solution snapshots at different time instances obtained from any black-box ODE solver in a library, at a small set of system parameter samples. In our approach, we fit a (simpler) time integration scheme of our choice to the available snapshots data. Since the time integration method we impose on the data does not coincide with the actual scheme used to generate the data in first place, there is a mismatch or \emph{defect} at each time step. Subject to a mild assumption, this defect corresponds to the local truncation error (LTE) of the imposed time integration scheme. We then learn this defect term as a function of time and the system parameters. Learning the defect using data allows us to formulate a new \emph{corrected} ROM that takes into account the LTE at each time step. Taking into account the LTE as a closure term means that the  corrected ROM can \emph{recover} the solution obtained from a ROM solved by a solver from an ODE library. We use the \emph{corrected} ROM to compute a good approximation to the true residual vector at each time step, leading to accurate estimation of the true error. The computation of the defect term is a key element of our proposed methodology. To do this efficiently, we rely on two observations. First, we demonstrate with numerical evidence that the defect term possesses a certain low-rank structure in space that allows us to efficiently project it on a low-dimensional subspace. Second, if the solution to the FOM is fairly smooth over the parameter space, then this smoothness carries over to the defect term and it can be efficiently approximated with respect to parameter variations using a suitable surrogate model.

\subsection{Prior work}
A number of recent works have sought to extend the RBM to a non-intrusive framework, mainly for steady systems; see \cite{morChaM09,morCasetal15,morChaH18}. In \cite{morChaH18} and \cite{morChaM09} the non-intrusive reduced basis (NIRB) method makes use of FEM solutions computed on two different meshes -- one fine and one coarse -- to estimate the state error due to the ROM. A different approach is proposed in \cite{morCasetal15} that involves the empirical interpolation method (EIM) to precompute an affine decomposition consisting of parameter-dependent and parameter-independent quantities. More recently, the method in \cite{morGroM22} has sought to extend the NIRB approach to time-dependent systems. All of these methods assume no knowledge of the FOM system matrices and rely purely on snapshots of the state variable. But, they need knowledge of the time integration scheme used to numerically integrate the spatially discretized system. In our new approach, we assume knowledge of the system matrices. Therefore, our approach is an intrusive method in terms of access to the model or the system matrices. However, our approach is non-intrusive in terms of the time integration scheme (ODE solver) used. Nevertheless, we do envisage an extension of our method to the case when there is no access to the system matrices; this would be a subject for future investigation. 

To eliminate the dependence on the time integration scheme, our approach aims to increase the accuracy of a user-imposed time integration method by addition of a data-driven closure term. The problem of improving time integration accuracy is an emerging topic and has received attention in the machine learning community~\cite{SheCL20,Polietal20,Huaetal22}. The deep Euler method (DEM) is introduced in \cite{SheCL20} where the authors' motivation is to improve the first-order accuracy of the explicit Euler time integration scheme by a factor of $\eta$, i.e., to $\mathcal{O}(\eta \delta t)$ ($\eta \ll 1$). They achieve this by approximating the LTE using a feed-forward neural network (FNN). Extension of the approach to other time integration schemes beyond the explicit Euler method are also illustrated. In the same spirit of the DEM approach, the work \cite{Polietal20} introduces \textsf{hypersolvers}, which are targeted at Neural ODEs~\cite{Cheetal18}. While both DEM and \textsf{hypersolvers} seek to improve the accuracy of the time integration scheme by approximating the LTE, they differ in the assumptions made about the model. DEM assumes the model of the PDE is known and is exact, while for \textsf{hypersolvers}, the model is unknown and is approximated by a separate neural network. Our work here differs from DEM and \textsf{hypersolvers} in several ways. Firstly, our aim is to obtain a corrected ROM for better error estimation whereas, in the aforementioned two methods the aim is simply to improve the FOM accuracy during simulation. Secondly, our targets are parametric nonlinear systems. Therefore, we want to learn the LTE for different time instances for a range of system parameters. Both DEM and \textsf{hypersolvers} are limited to non-parametric systems. While either method can potentially be extended to the parametric case, it is computationally very expensive. Furthermore, both these methods learn a closure term which is a function of the past states. Our approach treats the closure term as a function of time and any additional system parameters. This allows us to take advantage of the special structure and smoothness properties of the closure term.

The rest of this paper is organized as follows. In \Cref{sec:prelim}, we present the mathematical preliminaries of MOR and provide a short recap of the RBM along with \emph{a posteriori} error estimation. We also illustrate the pitfalls of the current error estimation approach used in RBM, through the example of the heat equation. \Cref{sec:improved_errestm} contains the main contributions of this manuscript. We introduce a user-imposed time integration scheme which incorporate the data-driven closure term and derive the \emph{a posteriori} output error estimator suited for situations where ODE solver libraries are used within the RBM. The algorithmic aspects of the proposed approach are discussed in \Cref{sec:computational_aspects} while numerical results are presented in \Cref{sec:numerics} to support our new method. We conclude with a summary and topics for future research in \Cref{sec:conclusion}.

\section{Mathematical background}%
\label{sec:prelim}

\begin{table}[t!]
	\small
	\caption{List of FOM variables}
	\centering
	\begin{center}
		\begin{tabular}{p{4cm}cc}
			\toprule
			\toprule
			Quantity & Variable & Dimension \\
			\toprule
			\toprule
			State vector & $\bx$ & $\R^{N}$ \\
			Input vector & $\bu$ & $\R^{N_{I}}$ \\
			Output vector & $\by$ & $\R^{N_{O}}$ \\
			Initial condition & $\bx_{0}$ & $\R^{N}$\\
			Nonlinear vector & $\bff$ & $\R^{N}$ \\		
			Approximate state vector & $\btx := \bV \bhx$ & $\R^{N}$ \\
			\midrule
			System matrix & $\bA$ & $\R^{N \times N}$ \\
			Input matrix & $\bB$ & $\R^{N \times N_{I}}$ \\
			Output matrix & $\bC$ & $\R^{N_{O} \times N}$ \\
			
			\bottomrule
		\end{tabular}
	\end{center}
	\label{tab:fom}
\end{table}
\begin{table}[t!]
	\small
	\caption{List of ROM variables}
	\centering
	\begin{center}
		\begin{tabular}{p{4cm}cc}
			\toprule
			\toprule
			Quantity & Variable & Dimension \\
			\toprule
			\toprule
			Left projection matrix & $\bV$ & $\R^{N \times n}$ \\
			Right projection matrix & $\bW$ & $\R^{N \times n}$ \\
			\midrule
			State vector & $\bhx$ & $\R^{n}$ \\
			Output vector & $\bhy$ & $\R^{N_{O}}$ \\			
			Initial condition & $\bhx_{0} := \bW^{\tpose} \bx$ & $\R^{n}$\\
			Nonlinear vector & $\bhf := \bW^{\tpose} \bff(\btx)$ & $\R^{n}$ \\
			\midrule
			System matrix & $\bhA := \bW^{\tpose} \bA \bV$ & $\R^{n \times n}$ \\
			Input matrix & $\bhB := \bW^{\tpose} \bB$ & $\R^{n \times N_{I}}$ \\
			Output matrix & $\bhC := \bC \bV$ & $\R^{N_{O} \times n}$ \\

			\bottomrule
		\end{tabular}
	\end{center}
	\label{tab:rom}
\end{table}

Consider the following parametric system of ODEs:
\begin{subequations}
	\label{eq:fom}
	\begin{align}
		\frac{d}{dt} \bx(\p) &= \bA(\p) \bx(\p) + \bff(\bx,\p) + \bB \bu(t), \qquad \bx(0) = \bx_{0},\label{eq:foma}\\
		\by(\p) &= \bC(\p) \bx(\p)\label{eq:fomb}.
	\end{align}
\end{subequations}
Such systems often arise upon discretizing a PDE using numerical discretization schemes. \Cref{tab:fom} lists each variable and its corresponding dimension. In subsequent discussions, we refer to \cref{eq:fom} as the FOM. The state vector $\bx(\p)$ can be obtained for any time instance for a given parameter $\p \in \mathcal{P} \subset \R^{p}$ by time integration of the FOM using any desired method like Runge-Kutta methods or linear multi-step methods. The set $\mathcal{P}$ denotes the parameter space and $p$ is its dimension.
The number of equations $N$ in \cref{eq:fom} is often very large to ensure a high-fidelity solution of the underlying physical process. This poses a major challenge for its numerical solution, especially for many instances of the parameter $\p$ in applications such as control and uncertainty quantification. 

\subsection{Model order reduction}
\label{subsec:mor}
Projection-based MOR techniques offer a systematic approach to obtain a ROM for \cref{eq:fom} in the following form:
\begin{subequations}
	\label{eq:rom}
	\begin{align}
		\frac{d}{dt} \bhx(\p) &= \bhA(\p) \bhx(\p) + \bhf(\btx, \p) + \bhB \bu(t), \qquad \bhx(0) = \bW^{\tpose} \bx_{0},\label{eq:roma}\\
		\bhy(\p) &= \bhC(\p) \bhx(\p)\label{eq:romb}.
	\end{align}
\end{subequations}
The above ROM is derived based on the ansatz $\bx \approx \btx := \bV \bhx$ applied to \cref{eq:fom}. The resulting over-determined system of equations is reduced by a Petrov-Galerkin projection, leading to \cref{eq:rom}. \Cref{tab:rom} lists the reduced variables and their corresponding dimensions. The number of equations $n$ in \cref{eq:roma} is often much smaller than $N$ in \cref{eq:foma}, i.e., $n \ll N$. Therefore, the ROM can be readily used for repeated simulations given any new values of the parameter $\p$. Different MOR techniques differ in how they compute the projection matrices $\bV, \bW$. 	When $\bV = \bW$, it is referred to as Galerkin projection. In the sequel, we will limit ourselves to ROMs obtained through a Galerkin projection.	

\begin{remark}
	Although \cref{eq:roma} is a ROM of dimension $n \ll N$, evaluating $\bhf(\btx, \p)$ involves operations that scale with the FOM dimension $N$ (since $\btx = \bV \bhx$ needs to be evaluated for each time step). Thus, a direct evaluation of \cref{eq:roma} may not offer any computational speedup over evaluating the FOM. Hyperreduction techniques~\cite{morBarMNetal04,morChaS10,morCarBF11} can be used to address this issue. We employ the discrete empirical interpolation method (DEIM) for the hyperreduction in our numerical results in \Cref{sec:numerics}.
\end{remark}

\subsection{Reduced basis method}
\label{subsec:rbm}

The RBM is a greedy approach that builds a global projection matrix $\bV$ to obtain a parametric ROM for \cref{eq:fom}. Since its introduction, the RBM has been a highly successful approach to obtain ROMs for parametric systems in a variety of applications such as process engineering~\cite{morZhaFLetal14,morCheFB19a}, geosciences~\cite{morDegVW20}, data assimilation~\cite{morDihH16,morKaretal18}, uncertainty quantification~\cite{morCheQR17} to name just a few. 

When RBM is applied to dynamical systems, the POD-Greedy method~\cite{morHaaO08} is adopted. It consists of a greedy sampling in the parameter space and a compression of the time trajectory through singular value decomposition (SVD). \Cref{alg:podgreedy} sketches the pseudo-code for the standard RBM using the POD-Greedy algorithm. In Step 1, \texttt{solver} denotes the time integration scheme used to solve \cref{eq:fom,eq:rom}. This is determined \emph{a priori} by the user based on the ODE library used, e.g. \texttt{scipy.ode.integrate}.

Consider a fine discretization of the parameter space $\mathcal{P}$ in the form of a training set $\Xi := \{\p_{1}, \p_{2}, \ldots, \p_{N_{p}}\}$ consisting of samples of the parameter $\p$. The POD-Greedy method starts by solving the FOM at a randomly chosen parameter $\p^{*} \in \Xi$. Performing the SVD of the solution snapshot matrix
\begin{align}
\label{eq:x_snapshot}
\bX(\p^{*}) := \big[ \bx(t_{0}), \bx(t_{1}), \ldots, \bx(t_{K})\big] \in \R^{N \times N_{t}}
\end{align}
with $N_{t} = K + 1$ yields the projection basis $\bV$. More precisely, we first obtain
\[
	\bX(\p^{*}) = \bV_{\p^{*}} \bm{\Sigma}_{\p^{*}} (\bW_{\p^{*}})^\tpose
\]
with $\bV_{\p^{*}} \in \R^{N \times N}, \bm{\Sigma}_{\p^{*}} \in \R^{N \times N_{t}}$ and $\bW_{\p^{*}} \in \R^{N_{t} \times N_{t}}$. The matrix $\bm{\Sigma}_{\p^{*}}$ is a rectangular diagonal matrix and contains the singular values $\sigma_{i}$ in the locations $\bm{\Sigma}_{\p^{*}}(i,\,i), i \in \{1, 2, \ldots, \min(N, N_{t})\}$ .

We then update $\bV$ by enriching it with the first $r_{c}$ left singular vectors of $\bX(\p^{*})$, i.e., $\bV_{\p^{*}}(: \,,\,1:r_{c})$. In subsequent iterations, snapshots of the FOM are similarly collected at different values of $\p^{*}$ and the matrix $\bV$ is updated with new information. Note, however, that for all iterations after the first, we perform a SVD of the matrix $\overline{\bX}$ obtained after removing from $\bX(\p^{*})$ the information already represented in $\bV$, i.e., we set
\[
	\overline{\bX} := 
	\bX(\p^{*}) - \bV (\bV^{\tpose} \bX(\p^{*})).
\]
To ensure good conditioning, it is recommended to perform a Gram-Schmidt orthonormalization after each update of $\bV$. The choice of $\p^{*}$ at each iteration is determined through an error estimator $\Delta(\p)$ as follows:
\[
	\p^{*} = \operatorname*{argmin}_{\p \in \Xi} \Delta(\p).
\]
The error estimator serves as an upper bound for the true state error $\| \bx(t, \p) - \btx(t, \p) \|$ (or true output error $\| \by(t, \p) - \bhy(t, \p) \|$). That is, 
\begin{align*}
	\| \bx(t, \p) - \btx(t, \p) \| &\leq \Delta(\p), \qquad \text{or}\\
	\| \by(t, \p) - \bhy(t, \p) \| &\leq \Delta(\p).
\end{align*}
All that is needed to evaluate $\Delta(\p)$ is to solve the ROM \cref{eq:roma} which can be formulated at each greedy iteration using $\bV$.
\begin{algorithm}[t]
	\small
	\caption{POD-Greedy algorithm}\label{alg:podgreedy}
	\begin{algorithmic}[1]
		\REQUIRE{Training set $\Xi$,
			tolerance ($\texttt{tol}$), 
			Discretized model ($\bE, \bA, \bB, \bC, \bff, \bx_{0}$)
		}
		
		\ENSURE{$\bV$
		}
		
		\vspace{0.5em}
		
		\STATE{Initialize $\bV = [\,\,]$,
			$\epsilon = 1 + \texttt{tol}$,
			choose \texttt{solver},
			greedy parameter $\p^{*}$
		}
		\vspace{0.5em}

		\WHILE{$\epsilon > \texttt{tol}$}
		\vspace{0.5em}
		\STATE{Obtain FOM snapshots $\bX(\p^{*})$ at $\p^{*}$ with \texttt{solver}}
		\vspace{0.5em}
		
		\STATE{Determine $\bV_{\p^{*}}$ through an SVD of $\overline{\bX} := 
			\bX(\p^{*}) - \bV (\bV^{\tpose} \bX(\p^{*}))$}
		\vspace{0.5em}
		
		\STATE Update $\bV$ as $\bV := \texttt{orth}\big( \bV, \bV_{\p^{*}}(: \,,\,1:r_{c})\big)$\\ with $\texttt{orth}\big(\cdot\big)$ denoting an orthogonalization process which can be implemented using the modified Gram-Schmidt process, or QR algorithm
		\vspace{0.5em}
		
		\STATE{Get ($\bhE, \bhA, \bhB, \bhC, \bhf, \bhx_{0}$) in \cref{eq:rom} by Galerkin proj. (+ hyperreduction)}
		\vspace{0.5em}
		
		\STATE{Solve ROM with \texttt{solver}}
		\vspace{0.5em}
		
		\STATE{$\p^{*} := \arg \max \limits_{\p \in \Xi} \Delta(\p)$}
		\vspace{0.2em}
		
		\STATE{Set $\epsilon = \Delta(\p^{*})$}
		\vspace{0.5em}
		
		\ENDWHILE
	\end{algorithmic}
\end{algorithm}
\subsection{A posteriori error estimation for the RBM}
\label{subsec:std_aposterr}
The standard error estimation approach in the RBM literature is residual-based~\cite{morQuaMN16, morHesRS16}.
In order to derive the error estimator $\Delta(\p)$, knowledge of the time discretization scheme used to integrate the FOM and the ROM is assumed, e.g., using implicit Euler, Crank-Nicolson method, or an implicit-explicit (IMEX) method. Computing $\Delta(\p)$ for a given parameter involves determining the residual vector $\br \in \R^{N}$ (or its norm) at each time instance; some approaches to error estimation also require the residual of a dual or adjoint system. Let us illustrate this by means of an example.

Suppose \cref{eq:fom} is discretized in time using a first-order IMEX scheme~\cite{AscRW95}. The linear part (involving $\bA(\p)$) is discretized implicitly, while the nonlinear vector $\bff(\bx)$ is evaluated explicitly. The resulting discretized system reads
\begin{subequations}
\label{eq:fom_timedisc}
\begin{align}
	\bE_{\text{im}} \bx_{\text{im}}^{k} &= \bA_{\text{im}} \bx_{\text{im}}^{k-1} + \delta t \big(\bff(\bx_{\text{im}}^{k-1}) + \bB \bu^{k} \big), \label{eq:fom_timedisca}\\
	\by_{\text{im}}^{k+1} &= \bC \bx_{\text{im}}^{k} \label{eq:fom_timediscb}
\end{align}
\end{subequations}
with $\bE_{\text{im}} := (\mathbf{I}_{N} - \delta t \bA)$, $\bA_{\text{im}} := \mathbf{I}_{N}$ where $\mathbf{I}_{N} \in \R^{N \times N}$ is the identity matrix.
For better clarity, we have not shown the parameter dependence of the system matrices and vectors. The ROM \cref{eq:rom} can be discretized in the same way and reads

\begin{subequations}
\label{eq:rom_timedisc}
	\begin{align}
	\bhE_{\text{im}} \bhx_{\text{im}}^{k} &= \bhA_{\text{im}} \bhx_{\text{im}}^{k-1} + \delta t \big(\bhf(\btx_{\text{im}}^{k-1}) + \bhB \bu^{k} \big), \label{eq:rom_timedisca}\\
	\bhy_{\text{im}}^{k+1} &= \bhC \bhx_{\text{im}}^{k} \label{eq:rom_timediscb}
	\end{align}
\end{subequations}
with $\bhE_{\text{im}} := (\mathbf{I}_{n} - \delta t \bhA)$, $\bhA_{\text{im}} := \mathbf{I}_{n}$ with $\mathbf{I}_{n} \in \R^{n \times n}$ being the identity matrix. The residual arising due to the ROM approximation can be computed by substituting the approximate state vector $\btx_{\text{im}}^{k}$ into \cref{eq:fom_timedisca}. The resulting residual at the $k$-th time step, $\br^{k}$, reads

\begin{equation}
\label{eq:resd_im1}
	\begin{aligned}
		\br^{k} := \bA_{\text{im}} \btx_{\text{im}}^{k-1} + \delta t \big(\bff(\btx_{\text{im}}^{k-1}) + \bB \bu^{k} \big) - \bE_{\text{im}} \btx_{\text{im}}^{k}.
	\end{aligned}
\end{equation}
It is clear that in order to obtain the residual $\br^{k}$, the time discretization scheme for the ROM should be the same as that for the FOM so that $\btx_{\text{im}}^{k}$ in \cref{eq:resd_im1} and $\bx_{\text{im}}^{k}$ in \cref{eq:fom_timedisc} correspond to the same time instance $t^{k}$. An \emph{a posteriori} error bound $\Delta(\p)$ for the approximation error $\be^{k}(\p) := \bx_{\text{im}}^{k}(\p) - \btx_{\text{im}}^{k}(\p)$ for a given parameter $\p$ can be computed based on the residual as below.
\begin{theorem}[Residual-based error bound]
	\label{thm:resd_err_bound_state}
	Suppose that the nonlinear quantity $\bff(\bx, \p)$ is Lipschitz continuous in the first argument for all $\p$ such that there exists a constant $L_{\bff}$ for which
	\[
	  \| \bff(\bx, \p) - \bff(\btx_{\textnormal{im}}, \p) \| \leq L_{\bff} \| \bx - \btx_{\textnormal{im}} \|.
	\]
	Further assume that for any parameter $\p$ the projection error at the first time step is $\| \be^{0}(\p) \| = \| \bx_{\textnormal{im}}^{0}(\p) - \btx_{\textnormal{im}}^{0}(\p) \| =\| \bx_{\textnormal{im}}^{0}(\p) - \bV \bV^{\tpose} \bx^{0}(\p) \|$. 
	
	The error in the state variable at the $k$-th time step, $\| \be^{k}(\p) \| = \| \bx^{k}_{\textnormal{im}}(\p) - \btx_{\textnormal{im}}^{k}(\p) \|$ is given by the following expression:
	\begin{align}
	\label{eq:std_errestm}
		\| \be^{k}(\p) \| \leq \Delta^{k}(\p) := \xi(\p)^{k} \nrm[2]{\mathbf{e}^{0}(\p)} + \sum\limits_{i=1}^{k} \zeta(\p) \cdot \big(\xi(\p)\big)^{k-i} \cdot \| \mathbf{r}^{i}(\p) \|_{2}
	\end{align}
	where $\zeta(\p) := \| \big(\bE(\p)\big)^{-1} \|_{2}$ and $\xi(\p) := 
	\Big( \| \big(\bE(\p)\big)^{-1} \bA(\p) \|_{2} + \delta t\,L_{\bff} \| \big(\bE(\p)\big)^{-1} \|_{2} \Big)$.
	\begin{proof}
		See \Cref{app:appendixA}.
	\end{proof}
\end{theorem}
Residual-based error bounds such as the one above are already available in the RBM literature~\cite{morHaaO08,morGre05} for both linear and nonlinear systems. 

In many applications, only a small set of variables which are obtained as a linear combination of the state variables are of interest. These are typically called the output quantities of interest. Goal-oriented \emph{a posteriori} error estimators for the output of interest have been discussed in several works~\cite{morZhaFLetal15,morCheFB19a,morGre05,morHaaO11}. They typically involve the residual of a dual or an adjoint system. The general form of the goal-oriented error estimators in~\cite{morGre05,morHaaO11} is
\begin{align}
	\label{eq:std_outerrestm}
	\| \by^{k}(\p) - \bhy^{k}(\p) \| \leq \Gamma(\p) \cdot \sum\limits_{i=1}^{k} \| \br^{i} \| \cdot \sum\limits_{j=k}^{K} \| \brdu^{j} \|
\end{align}
where $\| \brdu^{j} \|$ denotes the residual of an appropriately defined dual system at the $j$-th time instance and $\Gamma(\p)$ is some parameter-dependent constant. Improving upon this expression, a different goal-oriented error estimator was proposed in \cite{morZhaFLetal15,morCheFB19a} having the following form:
\begin{align}
\label{eq:zhang_outerrestm}
\| \by^{k}(\p) - \bhy^{k}(\p) \| \leq \breve{\Gamma}(\p) \cdot  \| \br^{k} \| \cdot  \| \brdu^{k} \|
\end{align}
with $\breve{\Gamma}(\p)$ being a parameter-dependent constant. This error estimator avoids the accumulation of the residuals over time, which is a major drawback for the goal-oriented error estimator in \cref{eq:std_outerrestm} and also the state error estimator in \cref{eq:std_errestm}.
It is also to be noted that the ROM resulting from the use of a goal-oriented error estimator often turns out to be of a smaller dimension.

Having reviewed the main ideas of RBM and the standard \emph{a posteriori} error estimator, we next illustrate the disadvantage of this standard approach when ODE solver libraries are used to solve the FOM/ROM within \Cref{alg:podgreedy}.

\subsection{RBM with ODE solver libraries}
Automatic ODE solver packages are implemented and readily available in a number of open-source and proprietary computational software. The \MATLAB $\,$ ODE Suite~\cite{matlabodesuite}, for instance, implements both linear multi-step (\texttt{ode15s, ode23s}, etc.) and Runge-Kutta-type solvers (\texttt{ode45, ode23}, etc.). In Python, using the \texttt{scipy} module once can access the \texttt{odeint} and \texttt{solve\_ivp} submodules, both providing access to a variety of standard time integration schemes. In addition to these, there are several stand-alone libraries to solve ODEs such as SUNDIALS~\cite{hindmarsh2005sundials}, ODEPACK~\cite{hindmarsh1983odepack}, ARKODE~\cite{reynolds2022arkode}. All of the mentioned libraries implement adaptive order and adaptive time-stepping which makes them highly efficient for a variety of problems, e.g., problems exhibiting stiffness.

When adaptive ODE solver packages are used within the RBM (i.e., in Steps 3 and 7 of \Cref{alg:podgreedy}, \texttt{solver} is replaced by the chosen method from the package), the standard error estimation approach (see \Cref{subsec:std_aposterr}) becomes less straightforward. One can no longer write the corresponding expression of the residual resulting from the ROM (e.g., \cref{eq:resd_im1}). The reason for this is that the exact expression of the time integration method used is unknown, as the \texttt{solver} adaptively varies the time step and/or the order of the scheme. If a residual expression obtained from a user-imposed time discretization scheme, e.g., \cref{eq:resd_im1} is used, the resulting error estimator \cref{eq:std_errestm} no longer gives an efficient bound for the true error, i.e., the error between $ \bV \bhx^{k}(\p)$ and $\bx^{k}(\p)$ obtained from simulating the ROM \cref{eq:rom} and the FOM \cref{eq:fom}, respectively, using an ODE solver from a library. We illustrate this next by obtaining a ROM for the simple non-parametrized linear heat diffusion model and estimating its error.

\subsection{Example: ROM for the linear heat equation}
\begin{figure}%
	\centering
	\includegraphics[scale=0.75]{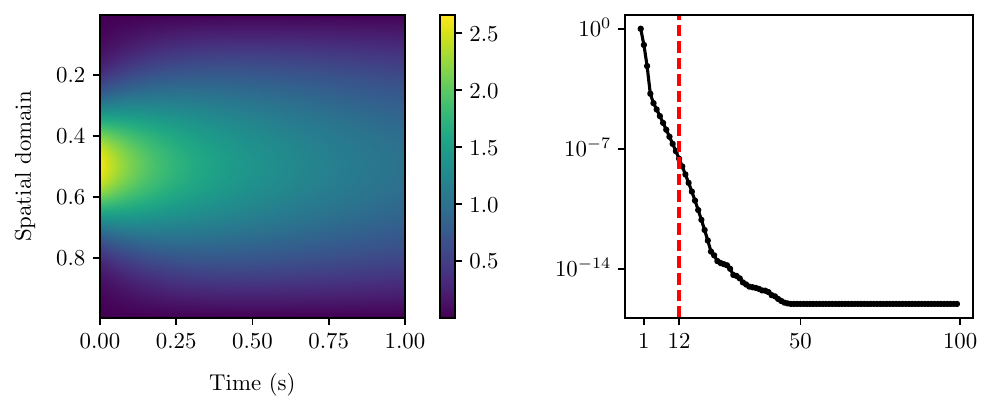}%
	\caption{Heat equation. Left: solution to the parametrized heat equation \cref{eq:fom_heat} at $\mu = 0.06$; Right: normalized singular values $\sigma_{i}/\sigma_{1}, i \in \{1, 2, \ldots, 100\}$.}
	\label{fig:heateqn_snapshot_svd}%
\end{figure}

We consider the linear heat equation in 1-D over the domain $\Omega = [0, 1]$ and time $t \in [0, 1]$
\begin{align}
	\label{eq:pde_heat_lin}
	\frac{\partial}{\partial t} v(z, t ; \mu) - \mu \frac{\partial^{2}}{\partial z^{2}}v(z, t ; \mu) = 0
\end{align}
where $v(z, t ; \mu)$ is the state, $z \in \Omega$ is the spatial variable and the viscosity $\mu = 0.06$. We further impose Dirichlet boundary conditions $v(0, t ; \mu) = v(1, t ; \mu) = 0,\, \forall t \in [0, 1]$. We also fix the output variable of interest as the value of the state at the node next to the right boundary. Employing the finite difference method, we discretize the domain $\Omega$ in an equidistant fashion, with a grid size of $2^{-8}$. With this, the resulting discretized ODE can be written in a form similar to \cref{eq:fom} as
\begin{align}
\label{eq:fom_heat}
\frac{d}{dt} \bx = \bA \bx, \qquad \bx(0) = \bx_{0}	
\end{align}
where $\bx \in \R^{255}$ is the discretized state vector and $\bA \in \R^{255 \times 255}$. We let the initial condition $\bx_{0}$ to be a normal distribution, i.e.,
\[
	\bx_{0} := \dfrac{1}{\sigma \sqrt{2 \pi}} e^{-\dfrac{1}{2} \cdot \bigg(\dfrac{z - \mathrm{m}}{\sigma}\bigg)^{2}}
\]
with a mean $\mathrm{m} = 0.5$ and a standard deviation of $\sigma = 0.15$. Our aim is to obtain a ROM for \cref{eq:fom_heat} and use the state and output error estimators in \cref{eq:std_errestm,eq:std_outerrestm}, respectively, to quantify the error of the ROM. For the output error estimator, we use the one proposed in \cite{morCheFB19a}. To this end, we consider the standard POD-based ROM. This involves obtaining snapshots of the state vector in \cref{eq:fom_heat} for different time instances. The SVD of the resulting snapshot matrix is used to obtain the projection matrix $\bV$. \Cref{fig:heateqn_snapshot_svd} illustrates the solution to \cref{eq:fom_heat} at $\mu = 0.06$. We also see the exponential decay of the singular values of the snapshot matrix. 

For this example, we compute $\bV$ consisting of the first $12$ columns of the left singular vector matrix. For integrating \cref{eq:fom_heat}, we use the \texttt{odeint} function available in the \texttt{scipy} package for Python. We note that \texttt{odeint} is a wrapper around the LSODA solver available in the Fortran library \texttt{odepack}. LSODA implements adaptive time-stepping. Moreover, it switches between methods for non-stiff and stiff problems automatically. We use the same solver to integrate the ROM. Once the reduced matrices are obtained through a Galerkin projection using $\bV$, the first step involved in estimating the error is to integrate the ROM using \texttt{odeint} to evaluate and plug-in the approximate solution $\widetilde{\bx}$ at each time instance, into the \texttt{odeint} time discretization scheme to get the residual. The residual operator at the $k$-th time instance has the general expression:
\begin{align}
	\mathcal{R}^{k}\big[\btx^{k}, \btx^{k-1}, \ldots, \btx^{k-s}\big]
\end{align}
where the arguments for $\mathcal{R}^{k}\big[\cdot\big]$ could be the solutions at the current and past $s$ time steps (in case the scheme used is a linear multi-step method) or at $s$ different stage solutions (in case an $s$-stage Runge-Kutta scheme is used). The exact form of this expression, naturally, is dependent on the time integration method that was used to obtain the snapshots. As we are using a solver from the \texttt{odepack} library, knowing the residual operator expression is complicated and often impossible. To circumvent this, one may choose to use a \emph{different, but known} time integration method in order to compute the residual, e.g., via \cref{eq:resd_im1}. However, this will lead to erroneous results, as we will demonstrate next. We denote by $\br^{k}$ the output of the residual operator $\mathcal{R}^{k}\big[\cdot\big]$ at a given set of arguments.
\begin{figure}%
	\centering
	\includegraphics[scale=0.75]{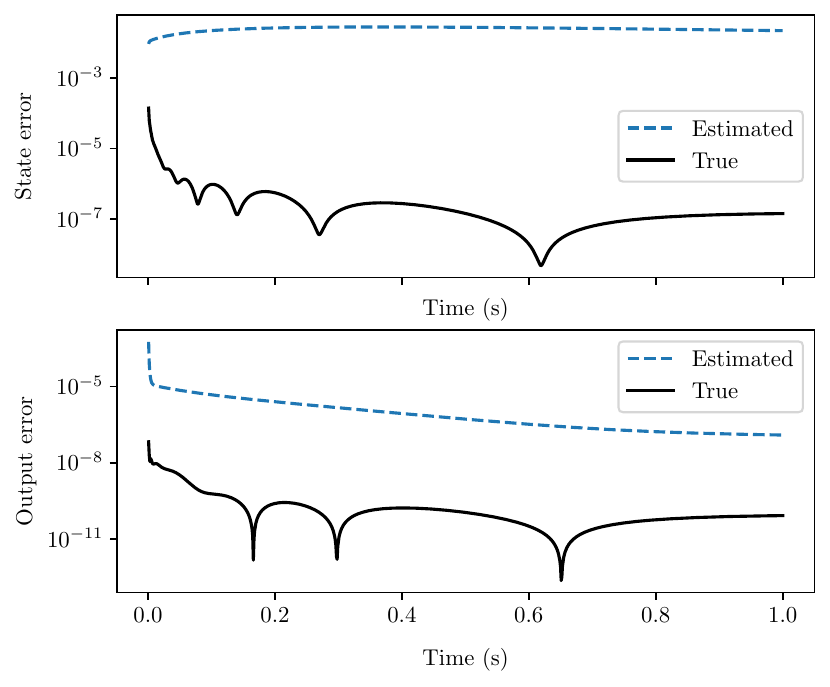}%
	\caption{Estimated and true errors for the heat equation \cref{eq:fom_heat} for $\mu = 0.06$ obtained by imposing a backward Euler method on the ROM snapshots obtained from \texttt{odeint}. Top: the estimated state error and the corresponding true error (see \cref{eq:std_errestm}); Bottom: the estimated output error and the corresponding true error (see \cref{eq:std_outerrestm}).}
	\label{fig:heateqn_outputerror}%
\end{figure}	

Suppose, we use a different time integration scheme, say a first-order backward Euler method as an ``approximation" to the true method used within \texttt{odeint}. In this case, the approximation to the true residual reads
\begin{align}
\label{eq:approx_residual_1storderEuler}
\widetilde{\br}^{k} =  \btx^{k-1} - \bigg( \mathbf{I}_{N} - \delta t \bA \bigg) \btx^{k}
\end{align}
with the approximate residual operator $\widetilde{\mathcal{R}}^{k}\big[\ast, \star] = \mathbf{I}_{N} \ast - \bigg( \mathbf{I}_{N}  - \delta t \bA \bigg) \star$ and $\ast, \star$ being placeholders for the arguments of $\widetilde{\mathcal{R}}$.
In general, $\mathcal{R}^{k}\big[\cdot\big]$ and $\widetilde{\mathcal{R}}^{k}\big[\cdot\big]$ are quite different, leading to rather inaccurate estimation of the error when the latter is used in \cref{eq:std_errestm}. This is illustrated in \Cref{fig:heateqn_outputerror}. The top figure shows the estimated state error obtained using \cref{eq:std_errestm} for every time step and the corresponding true error, both measured in the $2$-norm. The bottom figure illustrates the estimated error for the output variable and its true error.

The estimated error measured using the wrong expression of the residual overestimates by $2$ orders of magnitude in the best case. In this work, we propose a scheme to suitably modify $\widetilde{\mathcal{R}}^{k}\big[\cdot\big]$ using a closure term, such that the resulting expression for the residual is close to the one evaluated by $\mathcal{R}^{k}\big[\cdot\big]$. We limit our focus to the case of output error estimation though the proposed closure technique for correcting the residual can be straight-forwardly applied to any residual-based error estimator.

\section{Improving output error estimation via a data-enhanced closure approach}
\label{sec:improved_errestm}
In this section, we introduce a data-enhanced closure strategy to ensure that the residual resulting from the user-imposed time integration scheme, viz. $\widetilde{\mathcal{R}}^{k}\big[\btx^{k}, \btx^{k-1}, \ldots, \btx^{k-s}\big]$  is close to the true residual $\mathcal{R}^{k}\big[\btx^{k}, \btx^{k-1}, \ldots, \btx^{k-s}\big]$ such that the estimated output error is accurate.
\subsection{Defect-corrected FOM and ROM}
In our proposed approach, we first add a closure term to the FOM resulting from the user-imposed time integration scheme. This closure term is derived based on the snapshots of the true solution obtained using the ODE solver library. More precisely, suppose we have the snapshots of the solution \cref{eq:x_snapshot} at any given parameter $\p$ obtained from an ODE solver or some legacy codes:
\begin{align}
\label{eq:solver_snapshots}
\bX = \bigg[ \bx^{0} \, \, \bx^{1} \, \, \cdots \, \, \bx^{K}\bigg] \in \R^{N \times N_{t}}.
\end{align}
For purpose of illustration, we explain the details of the new method by considering a first-order IMEX scheme as the user-imposed time integration scheme. The FOM resulting from this is exactly \cref{eq:fom_timedisc}.

Since the user-imposed time integration scheme differs from the one used to generate the snapshots in $\bX$, we have a defect or a mismatch when we insert $\bx^{k}$ in \cref{eq:solver_snapshots} into the first-order IMEX scheme. This reads
\begin{align}
	\label{eq:defect_im1}
	\bd^{k} := \bE_{\text{im}} \bx^{k} - \Big(\bA_{\text{im}} \bx^{k-1} + \delta t \big( \bff(\bx^{k-1}) + \bB \bu^{k} \big) \Big).
\end{align}
\begin{remark}
	We note that this is precisely the local truncation error (LTE)~\cite{But08} of the first-order IMEX method if $\bx^{k}$ is assumed to be the true solution to \cref{eq:fom}. This assumption is not strong as the tolerances (absolute, relative) of the ODE solver can be set judiciously to achieve this.
\end{remark}

\begin{figure}[t!]
	\centering
	\includegraphics[scale=0.75]{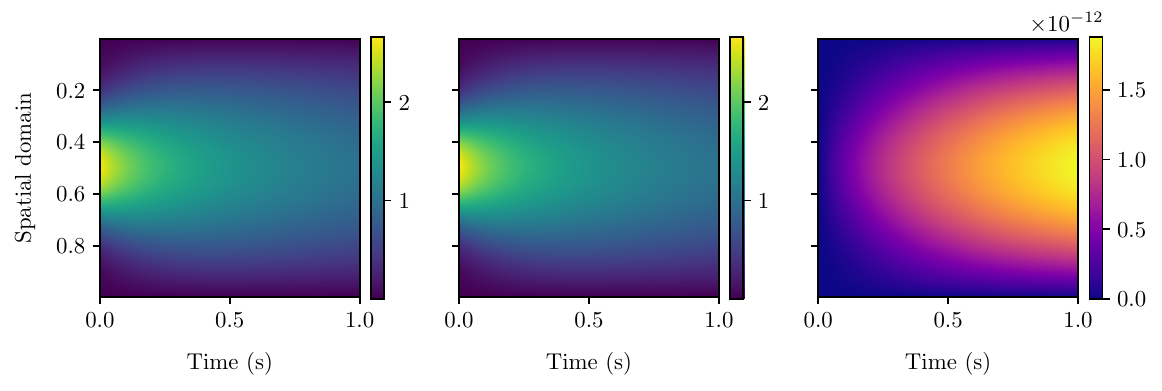}%
	\caption{Solution to the heat equation \cref{eq:fom_heat}; Left: solution  obtained using the ODE solver; Middle: solution using C-FOM \cref{eq:fom_timedisc_corr}; Right: pointwise errors between the solutions from the solver and C-FOM.}
	\label{fig:heateqn_FOm_CROM_compare}%
\end{figure}

We seek to modify the time-discrete FOM \cref{eq:fom_timedisc} such that its solution recovers the solution of \cref{eq:fom} computed by an ODE solver from a library. To this end, consider the following corrected FOM (C-FOM for short) obtained by adding the defect vector $\bd^{k}$ as a closure term:
\begin{subequations}
	\label{eq:fom_timedisc_corr}
	\begin{align}
		\bE_{\text{im}} \bx_{\text{im,c}}^{k} &= \bA_{\text{im}} \bx_{\text{im,c}}^{k-1} + \delta t \big(\bff(\bx_{\text{im,c}}^{k-1}) + \bB \bu^{k} \big) + \bd^{k}, \label{eq:fom_timedisca_corr}\\
		\by_{\text{im,c}}^{k} &= \bC \bx_{\text{im,c}}^{k}. \label{eq:fom_timediscb_corr}
	\end{align}
\end{subequations}
In \cref{eq:fom_timedisc_corr}, $\bx_{\text{im,c}}^{k}$ is the solution obtained after introducing the closure term and, as such, it differs from the solution $\bx_{\text{im}}^{k}$ to \cref{eq:fom_timedisc}. 
In fact, without adding $\bd^{k}$ as a closure term the local truncation error of the first-order IMEX method is $\mathcal{O}((\delta t)^{2})$. However, the local truncation error of the solution from the corrected IMEX method \cref{eq:fom_timedisc_corr} is of the same order as that resulting from the ODE solver. It can be seen from \Cref{fig:heateqn_FOm_CROM_compare} that the solution to the heat equation obtained using \cref{eq:fom_timedisc_corr} is identical to the one obtained using an ODE solver.

We emphasize that the defect vector $\mathbf{d}^{k}$ represents one choice for the closure term. In general, any other form of the closure term may be used. Since we use the defect vector as our choice for the closure term to recover the true residual, we use the two terms interchangeably.
\begin{remark}
	For any general time integration scheme the defect vector $\bd^{k}$ can be shown to have the following equivalence:
	\begin{align}
	\label{eq:defect_residual_equivalence}
	\bd^{k} = - \mathcal{R}^{k}\big[\bx^{k}, \bx^{k-1}, \ldots, \bx^{k-s}\big].
	\end{align}
\end{remark}

A corrected ROM (C-ROM) corresponding to the C-FOM can be defined by projecting the defect vector $\bd^{k}$ to the reduced space. The reduced defect vector is defined as $\widehat{\bd}^{k} := \bV^{\tpose} \bd^{k}$ and the C-ROM is:
\begin{subequations}
	\label{eq:rom_timedisc_corr}
	\begin{align}
	\bhE_{\text{im}} \bhx_{\text{im,c}}^{k} &= \bhA_{\text{im}} \bhx_{\text{im,c}}^{k-1} + \delta t \big(\bhf(\btx_{\text{im,c}}^{k-1}) + \bhB \bu^{k} \big) + \widehat{\bd}^{k}, \label{eq:rom_timedisca_corr}\\
	\bhy_{\text{im,c}}^{k} &= \bhC \bhx_{\text{im,c}}^{k}. \label{eq:rom_timediscb_corr}
	\end{align}
\end{subequations}

\subsection{An error estimator using the C-ROM}
Next, we make use of the C-ROM to derive a new error estimator $\overline{\Delta}^{k}(\p)$ that accurately estimates the true error $\| \by^{k}(\p) - \bhy^{k}(\p) \|$, where $\by^{k}(\p)$ and $\bhy^{k}(\p)$ are the output of the FOM \cref{eq:fom} and the ROM \cref{eq:rom} at time $t^{k}$, respectively. The FOM and and the ROM can be solved using any ODE solver. As mentioned, to derive the residual correctly, the same ODE solver must be applied to both the FOM and ROM simulations. Our proposed error estimator makes use of a dual system. We begin by deriving the dual system for \cref{eq:fom_timedisc_corr}.

\subsubsection{Dual system}
We derive the dual system corresponding to the C-FOM \cref{eq:fom_timedisc_corr} using the method of Lagrange multipliers. The Lagrangian can be formulated as
\begin{align}
	\mathcal{L} := \bC \bx_{\text{im,c}}^{k} + \big(\Lambda^{k}\big)^{\tpose} \bigg( \bE_{\text{im}} \bx_{\text{im,c}}^{k} - \bA_{\text{im}} \bx_{\text{im,c}}^{k-1} - \delta t \bff(\bx_{\text{im,c}}^{k-1}) - \delta t \bB \bu^{k} - \bd^{k} \bigg),
\end{align}
with $\Lambda \in \R^{N}$ being the vector of Lagrange coefficients. The dual system can be obtained by setting $\frac{\partial \mathcal{L}}{\partial \bx_{\text{im,c}}^{k}} \equiv \mathbf{0}$. This yields the system
\begin{align}
	\label{eq:fom_im1_dual}
	\bE_{\text{du}} \bx_{\text{du}} = \bC_{\text{du}}
\end{align}
with $\bE_{\text{du}} := \bE_{\text{im}}^{\tpose}$ and $\bC_{\text{du}} := -\bC^{\tpose}$. Note that the Lagrange multipliers $\Lambda$ are the dual state variables; for better clarity we denote them by $\bx_{\text{du}}$. Note that the defect vector $\mathbf{d}^{k}$ is treated as a function of time and parameter $\p$. Therefore, it does not depend on the solution $\bx_{\text{im,c}}^{k}$ of the C-ROM.

We further define the dual ROM as 
\begin{align}
	\label{eq:rom_im1_dual}
	\bhE_{\text{du}} \bhx_{\text{du}} = \bhC_{\text{du}}.
\end{align}
obtained by making the ansatz $\bx_{\text{du}} \approx \btx_{\text{du}} = \bV_{\text{du}} \bhx_{\text{du}}$. Here, $\bV_{\text{du}}$ is the projection matrix corresponding to the dual system and $\bhE_{\text{du}} := \bV_{\text{du}}^{\tpose} \bE_{\text{du}} \bV_{\text{du}}$, $\bhC_{\text{du}} := \bV_{\text{du}}^{\tpose} \bC_{\text{du}}$.

The residuals corresponding to the ROM \cref{eq:rom_timedisc_corr} and the dual ROM \cref{eq:rom_im1_dual} are, respectively,
\begin{align}
	\label{eq:resd_im1_corr}
	\br_{\text{im,c}}^{k} := \bA_{\text{im}} \btx_{\text{im,c}}^{k-1} + \delta t \big(\bff(\btx_{\text{im,c}}^{k-1}) + \bB \bu^{k} \big) + \bd^{k} - \bE_{\text{im}} \btx_{\text{im,c}}^{k}
\end{align}
and
\begin{align}
\label{eq:resddual_im1_corr}
\br_{\text{du}} := \bC_{\text{du}} - \bE_{\text{du}} \btx_{\text{du}}.
\end{align}
Following the approach in \cite{morZhaFLetal15,morCheFB19a} we define an auxiliary residual $\breve{\br}_{\text{im,c}}$ as
\begin{align}
	\label{eq:resdbr_im1_corr}
\breve{	\br}_{\text{im,c}}^{k} :&= \bA_{\text{im}} \bx_{\text{im,c}}^{k-1} + \delta t \big(\bff(\bx_{\text{im,c}}^{k-1}) + \bB \bu^{k} \big) + \bd^{k} - \bE_{\text{im}} \btx_{\text{im,c}}^{k}.
\end{align}
The auxiliary residual defined above will be required in the derivation of the data-enhanced output error estimator.

\subsubsection{Modified output term}
In anticipation of the new error estimator we propose next, we introduce a modified output variable $\overline{\by}_{\text{im,c}}^{k}$ defined as
\begin{align}
\label{eq:modified_output}
\overline{\by}_{\text{im,c}}^{k} := \bhy_{\text{im,c}}^{k} - \btx_{\text{du}}^{\tpose} \br_{\text{im,c}}^{k}.
\end{align}
Adding a correction term to the output in the form of a dual-weighted residual is an established practice~\cite{morRov03,morGre05,morCheFB19a} and serves to derive a tighter estimate of the true error.

\subsubsection{Data-enhanced error estimation}
Denoting the ODE solver applied to solve the FOM \cref{eq:fom} and the ROM \cref{eq:rom} as \texttt{solver}, the norm of the true error we desire viz., $\| \by^{k} - \bhy^{k} \|$ can be written as:
\begin{align}
	\| \by^{k} - \bhy^{k} \| &= \|   \by^{k} - \overline{\by}_{\text{im,c}}^{k} + \overline{\by}_{\text{im,c}}^{k} - \bhy^{k}  \| \\
	&\leq \| \by^{k} - \overline{\by}_{\text{im,c}}^{k} \| + \| \overline{\by}_{\text{im,c}}^{k} - \bhy^{k} \|. \label{eq:true_error_split}
\end{align}
The following theorem bounds the first summand in \cref{eq:true_error_split}:
\pagebreak
\begin{theorem}[A posteriori error bound for the corrected ROM]
	\label{thm:err_estm}
	Given the FOM in \cref{eq:fom}, the C-FOM in \cref{eq:fom_timedisc_corr} and the C-ROM \cref{eq:rom_timedisc_corr}, assuming that \emph{$\bE_{\text{im}}$} is non-singular for all $\p \in \mathcal{P}$, we have the following error bound for the modified output vector in \cref{eq:modified_output}:
	\begin{align}
		\| \by^{k} - \overline{\by}_{\textnormal{im,c}}^{k} \| \leq \| \bE_{\textnormal{im}}^{-1} \| \, \| \mathbf{r}_{\textnormal{du}} \| \, \| \breve{\mathbf{r}}_{\textnormal{im,c}}^{k} \| + \| \btx_{\textnormal{du}} \| \, \| \mathbf{r}_{\textnormal{im,c}}^{k} - \breve{\mathbf{r}}_{\textnormal{im,c}}^{k} \|.
	\end{align}
	
	\begin{proof}
		See \Cref{app:appendixB}.
	\end{proof}
\end{theorem}
Although the bound above is rigorous, it is not computable owing to the quantity $\breve{\br}^{k}$. Recall from \cref{eq:resdbr_im1_corr} that its computation requires that the FOM solution $\bx_{\text{im,c}}^{k}$ is available for any parameter $\p$, which is not the case. To derive a computable error estimator, we make use of the arguments used in \cite{morCheFB19a} to get the following error indicator:
\begin{align}
	\label{eq:err_estm_computable}
	\| \by^{k} - \overline{\by}_{\textnormal{im,c}} \| \lessapprox \bigg(\overline{\rho}\, \| \bE_{\textnormal{im}}^{-1} \|\, \| \mathbf{r}_{\textnormal{du}} \| + | 1 - \overline{\rho} | \, \| \btx_{\text{du}} \|\bigg) \, \| \br_{\text{im,c}}^{k} \|.
\end{align}
The quantity $\overline{\rho}$ is a measure for how close the residual $\br_{\text{im,c}}^{k}$ is to the auxiliary residual $\breve{\br}_{\text{im,c}}^{k}$ and is defined as
\begin{align}
	\label{eq:rhobar}
	\overline{\rho} = \frac{1}{K} \sum\limits_{k=1}^{K} \rho^{k}, \qquad \rho^{k} = \frac{\| \breve{\br}_{\text{im,c}}^{k}(\p^{*}) \|}{\| \br_{\text{im,c}}^{k}(\p^{*}) \|}.
\end{align}
It is evaluated only at the greedy parameter $\p^{*}$ for which snapshots of the true solution $\bx^{k}$ are available. For additional details, we refer to \cite{morCheFB19a}. Substituting \cref{eq:err_estm_computable} into \cref{eq:true_error_split} results in
\begin{align}
		\label{eq:data_enhanced_estimator_a}
		\| \by^{k} - \bhy^{k} \| &\lessapprox \bigg(\overline{\rho}\, \| \bE_{\textnormal{im}}^{-1} \|\, \| \mathbf{r}_{\textnormal{du}} \| + | 1 - \overline{\rho} | \, \| \btx_{\text{du}} \|\bigg) \, \| \br_{\text{im,c}}^{k} \| + \| \overline{\by}_{\text{im,c}}^{k} - \bhy^{k} \| =: \overline{\Delta}^{k}_{a}(\p).
\end{align}
The second quantity in the above inequality is the norm of the error between the output resulting from the user-imposed C-ROM and that obtained by solving the ROM \cref{eq:rom} with \texttt{solver}. It can be obtained cheaply as only two ROMs with small sizes need to be solved. Typically, our numerical experiments show that this quantity is very small and lesser than the magnitude of the first quantity. Therefore, we can safely neglect it so that only the first quantity can act as an alternative form of the proposed data-enhanced error estimator, i.e.,
\begin{align}
	\| \by^{k}(\p) - \bhy^{k}(\p) \| &\lessapprox \| \by^{k}(\p) - \overline{\by}_{\text{im,c}}^{k}(\p) \| \\ 
	&\lessapprox \bigg(\overline{\rho}\, \| \bE_{\textnormal{im}}^{-1} \|\, \| \mathbf{r}_{\textnormal{du}}(\p) \| + | 1 - \overline{\rho} | \, \| \btx_{\text{du}}(\p) \|\bigg) \, \| \br_{\text{im,c}}^{k}(\p) \| =: \overline{\Delta}^{k}_{b}(\p) \label{eq:data_enhanced_estimator}.
\end{align}
\subsubsection{Error estimation in presence of hyperreduction}
For nonlinear systems, the efficient computation of the ROM is impeded by the presence of the nonlinear function $\bff(\btx_{\text{im,c}}^{k})$ whose evaluation scales as the dimension $N$ of the FOM. To tackle this, the DEIM approach~\cite{morChaS10} is used in this work. Using DEIM, the residual expression in \cref{eq:resd_im1_corr} gets modified as below:
\begin{align}
	\br_{\text{im,c}}^{k} &= \bA_{\text{im}} \btx_{\text{im,c}}^{k-1} + \delta t \big(\bff(\btx_{\text{im,c}}^{k-1}) + \mathscr{I}[\bff(\btx_{\text{im,c}}^{k-1})] -  \mathscr{I}[\bff(\btx_{\text{im,c}}^{k-1})] + \bB \bu^{k} \big) + \bd^{k} - \bE_{\text{im}} \btx_{\text{im,c}}^{k},\nonumber\\
	&= \underbrace{\bigg(\bA_{\text{im}} \btx_{\text{im,c}}^{k-1} + \delta t \big( \mathscr{I}[\bff(\btx_{\text{im,c}}^{k-1})] +  \bB \bu^{k} \big) + \bd^{k} - \bE_{\text{im}} \btx_{\text{im,c}}^{k} \bigg)}_{\br_{\text{im,c}, \mathscr{I}}^{k}} +\, \underbrace{\delta t \bigg(\bff(\btx_{\text{im,c}}^{k-1}) - \mathscr{I}[\bff(\btx_{\text{im,c}}^{k-1})]\bigg)}_{\mathbf{e}_{H}^{k}}.
\end{align}
In the equation above, $\mathscr{I}[\cdot]$ denotes the hyperreduction operator. Further, $\mathbf{e}_{H}^{k}$ refers to the error introduced by hyperreduction at the $k$-th time step. For a detailed discussion on the computational aspects, we refer to the work~\cite{morCheFB19a}, where the simultaneous adaptive construction of the RBM and the DEIM bases vectors is also discussed. This adaptive bases construction approach is also implemented in our numerical experiments.

We make use of the new data-enhanced error estimators ($\overline{\Delta}^{k}_{a}(\p)$ or $\overline{\Delta}^{k}_{b}(\p)$) to choose the next parameter within the greedy algorithm. For each parameter $\p$, we determine the average estimated error over time as:
\begin{align}
	\label{eq:mean_estm_error}
	\overline{\Delta}_{z}(\p) = \dfrac{1}{N_{t}} \sum\limits_{k=0}^{N_{t}} \overline{\Delta}^{k}_{z}(\p),
\end{align}
where $z$ stands for $a$ or $b$. 
In the numerical results, we use $\overline{\Delta}^{k}_{b}(\p)$ in  \cref{eq:mean_estm_error} in Line 11 of  \Cref{alg:podgreedy_ode} so that during the greedy iterations, the ROM \cref{eq:rom} does not have to be repeatedly solved with \texttt{solver} in order to evaluate the second quantity in \cref{eq:data_enhanced_estimator_a}.

\begin{remark}
	While we have illustrated the data-enhanced error estimator using a first-order IMEX time integration scheme, any consistent time integration scheme can be used, including higher order ones. We will also demonstrate the use of a second-order IMEX scheme in the numerical examples. Of course, one should be cautious of the fact that a higher order time integration scheme comes with a larger computational cost.
\end{remark}

\begin{remark}
	When the defect vector $\mathbf{d}^{k}$ present in the residual (see \cref{eq:resd_im1_corr}) is known exactly, i.e., the C-ROM uses the same time integration method used for the FOM, we recover the \emph{a posteriori} output error estimator proposed in \cite{morCheFB19a}. This shows that our newly proposed error estimator is consistent with the case where the time integration method used is known.
\end{remark}

\begin{remark}
	Note that the C-ROM is not the finally derived ROM, but is only used within the greedy algorithm to derive the error estimator $\overline{\Delta}^{k}_{z}(\p)$ (see Steps 9-10 in \Cref{alg:podgreedy_ode}) which estimates the error between the solution of the FOM \cref{eq:fom} and that of the ROM \cref{eq:rom}, both being computed using any \texttt{solver}. The error estimator can also be used in the online stage where the user may wish to use any preferred \texttt{solver} to compute the ROM solutions. Moreover, given a good approximation of the defect vector, the solution of the C-ROM and the solution $\bhx^{k}$ from the ODE solver are nearly the same. Therefore, $\overline{\Delta}^{k}_{b}(\p)$ is almost as accurate as $\overline{\Delta}^{k}_{a}(\p)$.
\end{remark}

\section{Computational aspects}
\label{sec:computational_aspects}
While we have derived an accurate data-enhanced \emph{a posteriori} output error estimator, there remain a few computational challenges. In this section we highlight these challenges and propose efficient solutions to address them. 

In \cref{eq:data_enhanced_estimator} the residual term $\br_{\text{im,c}}^{k}$  involves the defect vector $\bd^{k}$ (see \cref{eq:resd_im1_corr}). However, determining this term involves knowing the true solution to \cref{eq:fom} obtained from the ODE solver at any parameter $\p$ (see \cref{eq:defect_im1}).
Without a cheap method to approximate $\bd^{k}$ for every parameter $\p$, the error estimator is not efficient. To alleviate this, we make several observations about the function $\bd(t, \p)$, which will lead to its efficient approximation. These observations relate to
\begin{enumerate}
	\item a certain low-rank structure (over space) that $\bd(t, \p)$ possesses and
	\item a smoothness over parameter variations that it inherits from the underlying parametric PDE.
\end{enumerate}

\subsection{Low-rank structure of the defect}
For a large class of problems, there exist a low-dimensional subspace onto which the solution snapshots can be projected, without incurring a large error. Indeed, this forms the underlying motivation for performing model order reduction using POD and other methods. Leveraging this fact, we assume that the defect trajectory at a given parameter $\p$, $\bd^{k}(\p), \forall k \in \{0, 1, \ldots, K\}$ can also be efficiently approximated by a suitable low-dimensional subspace.

\begin{figure}[t!]
	\centering
	\includegraphics[scale=0.75]{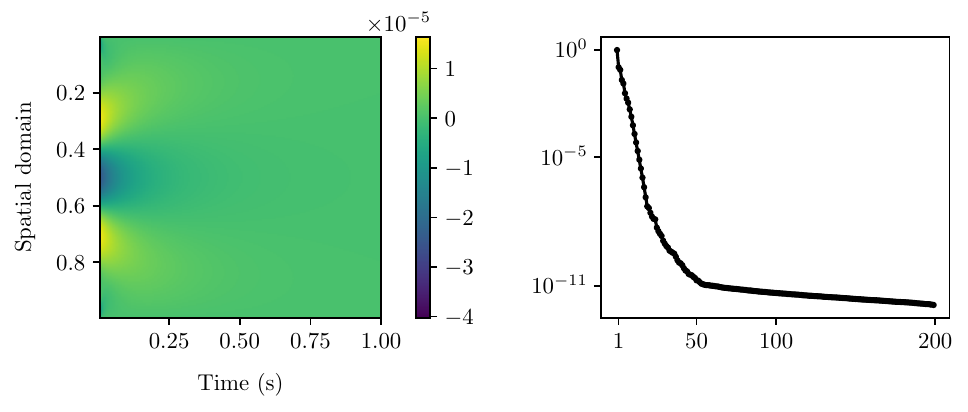}%
	\caption{Defect/local truncation error of the heat equation. Left: space time variation of the defect at $\mu = 0.06$; Right: singular values of the matrix $\mathbf{D} := \{\mathbf{d}^{k}\}_{k=0}^{K}$. The singular values exhibit an exponential decay illustrating the existence of a low-dimensional subspace.}
	\label{fig:heateqn_defect}%
\end{figure}

\Cref{fig:heateqn_defect} plots the values of the defect vector $\bd^{k}$ in \cref{eq:defect_im1} for the heat equation. It is evident that the defect trajectory has a certain regularity over the spatial domain (see left figure); evidence for the low-rank structure is also seen through the SVD performed (see right figure) on the snapshot matrix $\mathbf{D} \in \R^{N \times N_{t}}$ whose columns consist of snapshots of $\bd^{k}$ at different values of the time $t^{k}, \forall k \in \{0, \ldots, K\}$. We see an exponential decay of the relative quantity $\sigma_{i}/\sigma_{1}$, where $\sigma_{i}$ is the $i$-th singular value. A similar exercise will be repeated in \Cref{sec:numerics} for the three numerical examples we consider. In each of the cases, it will be seen that the defect trajectory is reducible spatially. We will use this fact to compute $\bd^{k}(\p)$ efficiently.

The concept of Kolmogorov $n$-width is used in the RBM literature~\cite{morQuaMN16,morBufetal12,Pin85} to quantify the approximability of the solution manifold $\mathcal{M}$ corresponding to a parametrized system of equations, with a linear subspace of dimension $n$ denoted as $\mathcal{V}^{n}$. Consider the solution manifold $\mathcal{M}$ for the FOM in \cref{eq:fom} defined as
\begin{align}
\label{eq:solution_manifold}
\mathcal{M} = \{ \bx(t, \p) \, : \, (t, \p) \in \mathcal{T} \times \mathcal{P} \} \subset \R^{N}.
\end{align}
The Kolmogorov $n$-width of $\mathcal{M}$ using $\mathcal{V}^{n}$ can be defined as
\begin{align}
	\label{eq:Knw_solution}
	d_{n}(\mathcal{M}) := \inf \limits_{\text{dim}(\mathcal{V}^{n})\,=\,n} \,\, \sup \limits_{\bx \in \mathcal{M}} \, \inf \limits_{\btx \in \mathcal{V}^{n}} \| \bx - \btx \|.
\end{align}
For the parametrized defect function $\bd(t, \p)$, we define the following manifold:
\begin{align}
	\label{eq:defect_manifold}
	\mathcal{M}_{d} = \{ \bd(t, \p) \, : \, (t, \p) \in \mathcal{T} \times \mathcal{P} \} \subset \R^{N}.
\end{align}
From our numerical examples we have observed that $\mathcal{M}_{d}$ can be well-approximated by low-dimensional subspaces when the original parametric problem has fast Kolmogorov $n$-width decay. This implicates that the Kolmogorov $n$-width of $\mathcal{M}_{d}$ may inherit the behaviour of $d_{n}(\mathcal{M})$. At present we can not strictly prove this, but from the definition of the defect vector in \cref{eq:defect_im1}, the defect vector $\bd(t, \p)$ can be seen as being in the image of the operator $\mathcal{D} : \mathbb{R}^{n} \times \mathcal{P} \mapsto \R^{n}$ and $\mathcal{D}[\bx(t, \p)] = \bE_{\text{im}} \bx(t, \p) - \big( \bA_{\text{im}} \bx(t-\delta t, \p) + \delta t (\bff(\bx(t-\delta t, \p)) + \bB \bu(t)) \big)$. With this observation, the inherit property of the Kolmogorov $n$-width of $\mathcal{M}_{d}$ might be proved based on Theorem 4.1 in \cite{CohD15}. We leave this as our future work.

\subsection{Strategies to approximate the defect}
As discussed in the previous section, the defect vector is assumed to possess a certain low-rank structure in space and it inherits the smoothness of the solution $\bx(t, \p)$ over parameter variations, from the underlying parametric PDE. We use these two observations to make the approximation of $\bd^{k}(\p)$ computationally efficient, such that the output error estimator defined in \cref{eq:data_enhanced_estimator} can be used in the POD-Greedy algorithm. To this end, we adopt a two-stage approximation strategy. 

Starting from the observation about the low-rank structure of the defect, we can approximate the defect vector at a given time and parameter using the basis expansion
\begin{align}
\label{eq:defect_basis_expansion}
	\bd(t,\p) \approx \widetilde{\bd}(t,\p) = \sum_{i=1}^{n_{d}} \mathbf{v}_{d, i}\, \widehat{d}_{i}(t,\p)
\end{align} 
where $\bv_{d,i} \in \R^{N}$ are the expansion bases and $\widehat{\bd}(t,\p) := [\widehat{d}_{1}(t, \p), \ldots, \widehat{d}_{n_{d}}(t, \p)] \in \R^{n_{d}}$ is the vector of expansion coefficients. We denote with $\bV_{d} := \big[ \bv_{d,1}, \bv_{d,2}, \ldots, \bv_{d,n_{d}} \big] \in \R^{N \times n_{d}}$ the basis matrix. If the observation regarding the rapid decay of the singular values holds, then $n_{d}$ will be small.
Given such a basis expansion for the defect vector, we can approximate it for any given parameter $\p$ and a time instance $t$ if $\widehat{\bd}(t,\p)$ can be evaluated cheaply. Our two-stage approach involves:
\begin{itemize}
	\item identifying a suitable basis matrix $\bV_{d}$ using a POD/SVD-based approach and
	\item learning the map $(t, \p) \mapsto \widehat{\bd}(t, \p)$ for which we propose two different approaches: one based on radial basis function interpolation and the other using a feed-forward neural network.
\end{itemize}

\subsubsection{SVD-based spatial reduction}
\label{subsec:SVD_defect_reduction}
In the first stage of approximation, we collect snapshots of the defect vector $\bd(t, \p), \forall t \in \mathcal{T}_{d} := \{t^{1}, t^{2}, \ldots, t^{N_{t}}\}$ and $\p \in \Xi_{\text{defect}}$ where $\Xi_{\text{defect}}$ is a set containing $d_{s}$ parameter samples, with $d_{s}$ typically small. Doing so involves solving \cref{eq:fom} with \texttt{solver} to obtain the FOM solution snapshots. Following this, the solution snapshots can be used to obtain the defect vector from e.g., \cref{eq:defect_im1}\footnote{Note that the defect is obtained here for a particular user-imposed time integration scheme.}. We denote by $\mathfrak{D} \in \R^{N \times N_{t} \times d_{s}}$ the third-order tensor arranged such that each of its frontal slices corresponds to the matrix $\mathbf{D}(\p) := \big[\mathbf{d}(t^{1}, \p), \mathbf{d}(t^{2}, \p), \ldots, \mathbf{d}(t^{N_{t}}, \p) \big] \in \R^{N \times N_{t}}$ and $\p \in \Xi_{\text{defect}}$. We refer to the $i$-th frontal slice as $\mathfrak{D}^{(i)}$. Next, we apply a two-step SVD reduction~\cite{morWanHR19} to $\mathfrak{D}$ which will result in $\widehat{\mathfrak{D}} \in \R^{n_{d} \times N_{t} \times d_{s}}$ as follows:
\begin{itemize}
	\item In step 1, we perform SVDs for each frontal slice separately, i.e., we perform a SVD for each $\mathfrak{D}^{(i)}$ with $i \in \{1, 2, \ldots, d_{s}\}$. For every such SVD, we make use of a fixed tolerance $\texttt{tol}_{\text{SVD}, t}$ to truncate the singular values and collect the first $\ell_{i}$ left singular vectors in the matrix $\bU_{d, i} \in \R^{N \times \ell_{i}}$.
	
	\item In step 2, we form a matrix $\mathbf{R}$ defined as 
	\[
		\mathbf{R} := \bigg[ \bU_{d, 1} \, \big\lvert \, \bU_{d, 2} \, \big\lvert \, \cdots \, \big\lvert \, \bU_{d, d_{s}}\bigg] \in \R^{N \times (\ell_{1} + \cdots + \ell_{d_{s}})}
	\]
	whose columns consist of the truncated left singular vectors for each parameter obtained in step 1. We then perform the SVD of $\mathbf{R}$, using a tolerance $\texttt{tol}_{\text{SVD}, \p}$ to obtain the projection matrix $\bV_{d} \in \R^{N \times n_{d}}$ as the first $n_{d}$ columns of the left singular vectors of $\mathbf{R}$. 
\end{itemize}
Finally, the reduced tensor $\widehat{\mathfrak{D}}$ can be obtained via a \textit{mode-1} tensor-matrix product as
\[
	\widehat{\mathfrak{D}} = \mathfrak{D} \times_{1} \bV_{d}^{\tpose}.
\]
We note that each mode-1 fiber in $\widehat{\mathfrak{D}}$ corresponds to the reduced defect vector $\widehat{\bd}(t^{k}, \p) \in \R^{n_{d}}$.

Thus far, we have reduced the dimension of the first mode of the tensor $\mathfrak{D}$ from $N$ to $n_{d}$, with $n_{d} \ll N$. Next, we detail the two approaches used to approximate $\widehat{\bd}(t, \p)$.
\begin{remark}
	The motivation for using the two-step SVD approach is to reduce the overall computational costs. While we have not pursued it in this work, another valid approach to reduce the computational cost would be a randomized-SVD.
\end{remark}

\subsubsection{Interpolation using radial basis functions}
Radial basis functions (RBFs) are a popular class of kernel methods which are used in scattered-data approximation~\cite{Buh03,Wed05}. We use RBFs to learn an approximation $\widehat{d}_{\text{RBF}}(\p)$ to each entry of the reduced defect vector at each time instance, such that for any given $(j, k)$ pair, $\widehat{d}_{\text{RBF}}(\p)$ interpolates $\widehat{d}_{j}(t^{k}, \p)$, the $j$-th entry of the reduced defect $\widehat{\bd}(t^{k}, \p)$ at parameters $\p \in \Xi_{\text{defect}}$. 

The RBF approximation reads
\begin{align}
\label{eq:rbf_basis_expansion}
	\widehat{d}_{\text{RBF}, j}^{k}(\p) = \sum\limits_{i=1}^{d_{s}} w_{i} \Phi(\| \p - \p_{i}  \|)
\end{align}
with $\{w_{i}\}_{i=1}^{d_{s}}$ denoting the weights and $\Phi(\cdot)$ being the radial basis functions. The weights are obtained by imposing the interpolation condition $\widehat{d}_{\text{RBF},j}^{k}(\p_{i}) = \widehat{d}_{j}(t^{k}, \p_{i}), \, i = 1, 2, \dots, d_{s}$, $\p_{i} \in \Xi_{\text{defect}}$. This leads to the following system of linear equations
\begin{align}
\label{eq:rbf_linear_syste}
		\begin{bmatrix}
		\Phi(\nrm{\p_{1}- \p_{1}}) & \Phi(\nrm{\p_{1}- \p_{2}}) &\cdots &\Phi(\nrm{\p_{1}- \p_{d_{s}}})\\
		\Phi(\nrm{\p_{2} - \p_{1}}) & \Phi(\nrm{\p_{2}- \p_{2}}) & \cdots & \Phi(\nrm{\p_{2}- \p_{d_{s}}})\\
		\vdots & \vdots & \ddots & \vdots\\
		\Phi(\nrm{\p_{d_{s}}- \p_{1}}) & \Phi(\nrm{\p_{d_{s}}- \p_{2}}) & \cdots & \Phi(\nrm{\p_{d_{s}}- \p_{d_{s}}})
		\end{bmatrix}
	\begin{bmatrix}
		w_{1} \\ 
		w_{2} \\
		\vdots\\
		w_{d_{s}}
		\end{bmatrix} =
	\begin{bmatrix}
		\widehat{d}_{j}(t^{k}, \p_{1})\\
		\widehat{d}_{j}(t^{k}, \p_{2})\\
		\vdots\\
		\widehat{d}_{j}(t^{k}, \p_{d_{s}})
		\end{bmatrix}.
\end{align}

\begin{figure}
	\centering
	\includegraphics[scale=0.35]{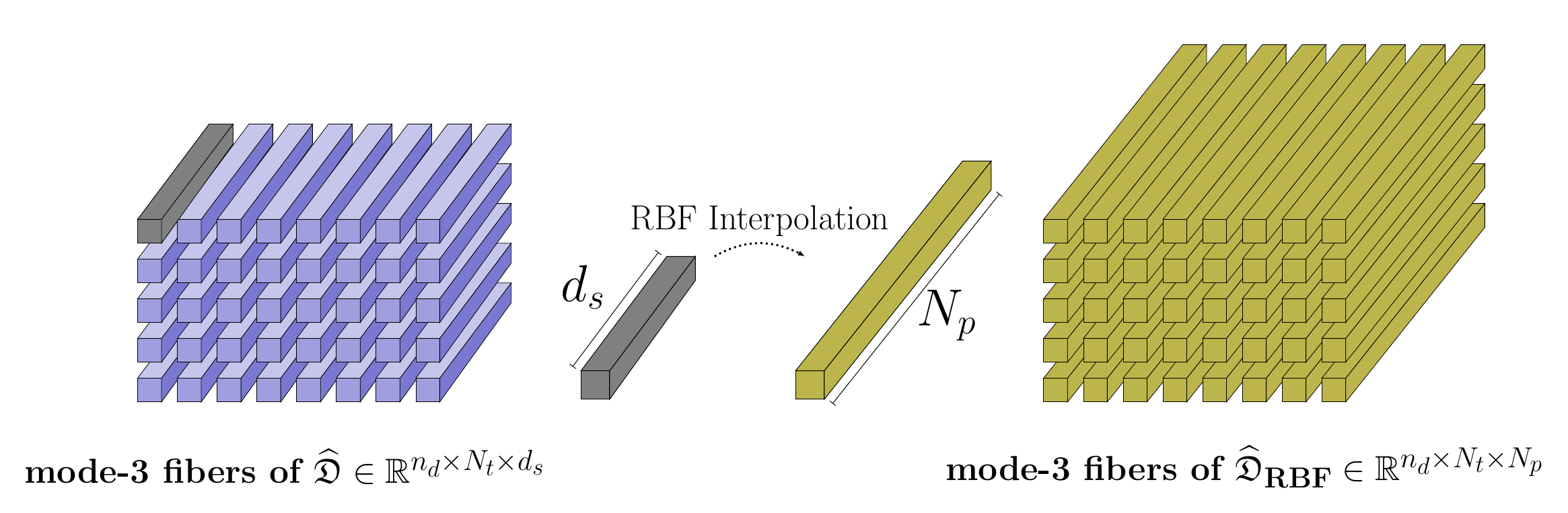}
	\caption{Approximation of the defect as a function of the parameter $\p$. The RBF interpolant learns an approximation of the defect vector over $N_{p}$ samples, with interpolation occurring at $d_{s}$ samples. We construct an individual RBF interpolant for each time and generalized spatial coordinate.}
	\label{fig:defect_approx_RBF}
\end{figure}

Based on the RBF interpolation, the defect vector $\bd(t^{k}, \p)$ in \cref{eq:defect_basis_expansion} can be approximated as
\begin{align}
	\label{eq:defect_RBF_approx}
	\bd(t^{k}, \p) \approx \widetilde{\bd}_{\text{RBF}}(t^{k}, \p) := \sum_{j=1}^{n_{d}} \mathbf{v}_{d, j}\, \widehat{d}^{k}_{\text{RBF},j}(\p)
\end{align}
We obtain an RBF interpolant for each time instance $k \in \{0, 1, 2, \ldots, K\}$ and each coordinate $j \in \{1, 2, \ldots, n_{d}\}$, resulting in a total of $(n_{d} \cdot N_{t})$ RBF interpolants. \Cref{fig:defect_approx_RBF} graphically illustrates the approach. Theoretically, $\widetilde{\bd}_{\text{RBF}}(t^{k}, \p)$ in \cref{eq:defect_RBF_approx} is valid for any $\p \in \mathcal{P}$. Therefore, $N_{p}$ in \Cref{fig:defect_approx_RBF} can be arbitrarily large. We denote it by $N_{p}$ as we are only interested in obtaining the reduced defect coefficients corresponding to the $N_{p}$ parameter samples present in the training set $\Xi$.

In our numerical results in \Cref{sec:numerics}, we denote this method of approximating the defect vector using a SVD spatial reduction followed by an RBF interpolation as \textsf{SVD+RBF}.

\begin{remark}
	In this work, we have considered separate RBF interpolants for each time step and each generalized spatial coordinate, leading to potentially many interpolants. While the number of interpolants scales as $\mathcal{O}(n_{d} N_{t})$, this can be efficiently implemented in one step by solving a linear system with the small $d_{s} \times d_{s}$ coefficient matrix in \cref{eq:rbf_linear_syste}, and with multiple right-hand sides (totalling $n_{d} N_{t}$). The runtime of solving such a linear system is usually much faster than separately solving $n_{d} N_{t}$ linear systems with the same coefficient matrix.
\end{remark}
\subsubsection{Approximation using artificial neural networks}
\begin{figure}[t!]
	\centering
	\includegraphics[scale=0.35]{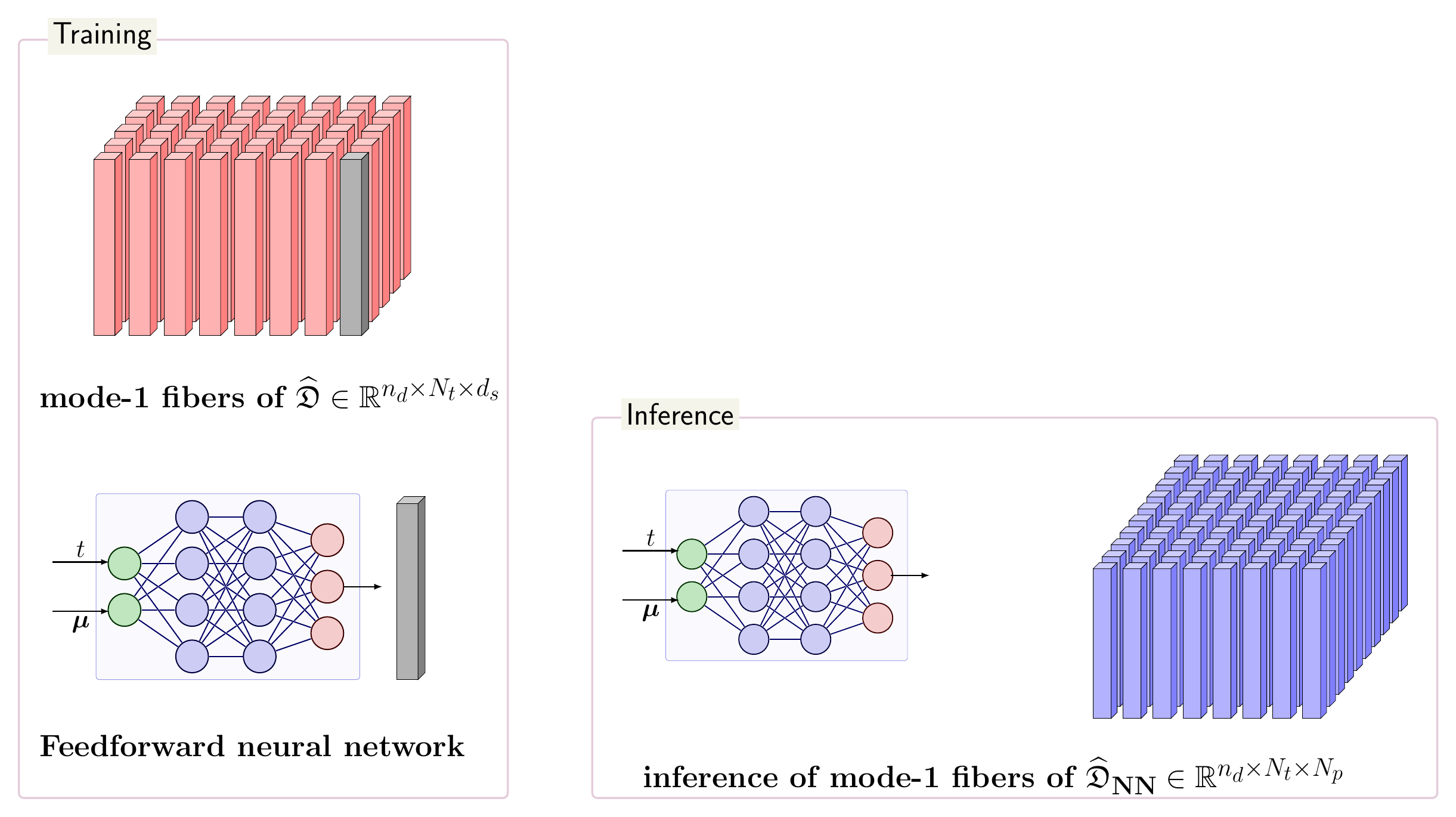}
	\caption{Approximation of the defect coefficients as a function of the inputs $(t, \p)$. The neural network is trained based on data available at $d_{s}$ parameter samples. In the inference stage, the neural network learns the approximation of the defect for all $N_{p}$ parameter samples.}
	\label{fig:defect_approx_NN}
\end{figure}
The second approach we consider to approximate the expansion coefficients $\widehat{\bd}(t, \p)$ for different time and parameter values is based on artificial neural networks (ANNs). A widely used architecture to implement ANNs are the feed-forward neural networks (FNNs). FNNs have been shown to be efficient for both regression and classification tasks in a variety of applications. The basic architecture of a FNN consists of three components: an input layer, hidden layer(s) and an output layer. The core element of any NN in general and FNNs in particular are artificial neurons. The hidden layer(s) in a FNN consists of neurons stacked together. Each neuron can receive inputs from a previous layer. The output corresponds to a nonlinear function of the weighted sum of its input signals. For a detailed overview of ANNs and FNNs, we refer to \cite{Bis06}. 

In this work, we consider the FNN to learn an approximation to the map between the inputs $(t, \p)$ and the output $\widehat{\bd}(t, \p) \in \R^{n_{d}}$. That is,
\[
\widehat{\bd}_{\text{NN}} : (t, \p) \mapsto \widehat{\bd}_{\text{NN}}(t, \p) \approx \widehat{\bd}(t, \p).
\]
To train the FNN, our training data consists of the dataset $\tau_{\text{train}} := (\bm{\Lambda}_{\text{NN}}, D_{\text{NN}})$ where $\bm{\Lambda}_{\text{NN}} := \{(t_{i}, \p_{j})\}_{i=j=1}^{N_{t}, d_{s}}$ is the input set, containing time and parameter samples and $D_{\text{NN}} := \{\widehat{\bd}(t_{i}, \p_{j})\}_{i=j=1}^{N_{t}, d_{s}}$ is the output data set which consists of the reduced defect $\widehat{\bd}(t_{i}, \p_{j}) \in \R^{n_{d}}$ at each time and parameter sample of the input. The neural network is implemented in PyTorch; more details regarding the number of layers used and other hyperparameters will be provided in the numerical section. The loss function is the mean square loss, viz.,
\[
	L_{\text{mse}} := \frac{1}{2} \sum_{i=1}^{d_{s}} \sum_{k=1}^{N_{t}} \bigg \lvert \bigg \lvert \widehat{\bd}(t^{k}, \p_{i}) - \widehat{\bd}_{\text{NN}}(t^{k}, \p_{i}) \bigg \rvert \bigg \rvert_{2}^{2}.
\] 
Once the neural network is trained, it can infer the values of the defect vector $\widehat{\bd}_{\text{NN}}(t, \p)$ at any chosen $(t, \p)$. The original defect vector at a given parameter $\p$ and at any time instance $\mathbf{d}(t, \p)$ can be approximated using the FNN-based approach as
\begin{align}
	\label{eq:defect_NN_approx}
	\bd(t, \p) \approx \widetilde{\bd}_{\text{NN}}(t, \p) := \bV_{d}\, \widehat{\bd}_{\text{NN}}(t, \p).
\end{align}
In our numerical results in \Cref{sec:numerics}, we denote this method of approximating the defect vector first with a SVD spatial reduction followed by an approximation of the coefficients $\widehat{\mathbf{d}}$ with a feedforward neural network as \textsf{SVD+FNN}.
\subsection{POD-Greedy with black-box ODE solvers}
We present the POD-Greedy algorithm that supports black-box ODE solvers and  incorporates the new data-enhanced \emph{a posteriori} error estimator as \Cref{alg:podgreedy_ode}. It requires some additional inputs compared to \Cref{alg:podgreedy}. These include a separate training set $\Xi_{\text{defect}} \subset \Xi$ for the defect approximation and two separate tolerances for the SVD, $\texttt{tol}_{\text{SVD}, t}, \texttt{tol}_{\text{SVD}, \p}$ which are required for the two-step SVD method to compute $\bV_{d}$. Before starting of the greedy algorithm in Step 4, Steps 2-3 in \Cref{alg:podgreedy_ode} are targeted towards learning the defect vector, which is added as a closure term to get the C-ROM. In Step 2, FOM solutions of \cref{eq:fom} are obtained for the $d_{s}$ parameter samples in the training set $\Xi_{\text{defect}}$. Using this data, the defect vectors (see \cref{eq:defect_im1}) induced by a user-imposed time integration scheme (denoted $\texttt{solver}_{\text{imp}}$) are collected in the tensor $\mathfrak{D}$. Then, in Step 3, the data in $\mathfrak{D}$ are first compressed into a low-dimensional space and the reduced defect vectors $\widehat{\mathbf{d}}(t, \p)$ for all $\p \in \Xi$ are learned via RBF or FNN. Afterwards, the defect vector $\mathbf{d}(t, \p)$ is approximated via the decoded vectors $\widetilde{\bd}_{\text{RBF}}(t, \p)$ or $\widetilde{\bd}_{\text{NN}}(t, \p)$. The approximation $\widetilde{\bd}_{\text{RBF}}(t, \p)$ or $\widetilde{\bd}_{\text{NN}}(t, \p)$ can be updated (replaced) by the true defect vector once new FOM data is available at $\p^{*}$ (Step 6, Step 11). Updating the approximate defect vector with the available true defect vector at $\p^{*}$ at each iteration leads to considerable improvements in the performance of \Cref{alg:podgreedy_ode}, as we shall demonstrate in \Cref{sec:numerics}. Furthermore, it reduces the amount of initial training samples needed in $\Xi_{\text{defect}}$. Typically, the user may not know, \textit{a priori}, the number of FOM samples needed to get a good approximation of the closure term. Therefore, $\Xi_{\text{defect}}$ can be coarsely sampled to keep the computational cost low. Once the greedy algorithm is begun, the FOM solution snapshots at $\p^{*}$ chosen at each greedy iteration are readily available. Those snapshots can be further used to compute the true defect vector $\mathbf{d}(t^{k}, \p^{*})$ via \cref{eq:defect_im1}. Since the snapshots are anyway available, the only computational costs incurred are those corresponding to the evaluation of $\mathbf{d}(t^{k}, \p), \forall k \in \{0, 1, \ldots, K\}$.
\begin{algorithm}[t]
	\small
	\caption{POD-Greedy algorithm for ODE solver libraries\label{alg:podgreedy_ode}}
	\begin{algorithmic}[1]
		\REQUIRE{Training set $\Xi$,
			tolerance ($\texttt{tol}$), 
			Discretized model ($\bE, \bA, \bB, \bC, \bff, \bx_{0}$)
		}
		\renewcommand{\algorithmicrequire}{\textbf{Input for defect approximation:}}
		\vspace{0.5em}
		\REQUIRE{Training set for defect approximation $\Xi_{\text{defect}} \subset \Xi$,
			SVD tolerances ($\texttt{tol}_{\text{SVD}, t}, \texttt{tol}_{\text{SVD}, \p}$),
			user-imposed time integration method $\texttt{solver}_{\textnormal{imp}}$
		}
		\vspace{0.5em}
		\ENSURE{$\bV$.
		}
		
		\vspace{0.5em}
		
		\STATE{Initialize $\bV = [\,\,]$,
			$\epsilon = 1 + \texttt{tol}$,
			choose \texttt{solver} from a library,
			greedy parameter $\p^{*}$
		}
		\vspace{0.5em}
		
		\STATE{Compute FOM snapshots for all $\p \in \Xi_{\text{defect}}$; determine the defect data tensor $\mathfrak{D}$ (see \Cref{subsec:SVD_defect_reduction})}
		\vspace{0.2em}
		\STATE{Perform SVD on $\mathfrak{D}$ to get $\bV_{d}$; obtain an approximation $\widetilde{\bd}_{\text{RBF}}(t, \p)$ \cref{eq:defect_RBF_approx} or $\widetilde{\mathbf{d}}_{\text{NN}}(t, \p)$ \cref{eq:defect_NN_approx} for $\mathbf{d}(t, \p)$}
		
		\WHILE{$\epsilon > \texttt{tol}$}
		\vspace{0.5em}
		\STATE{Obtain FOM snapshots $\bX_{\p^{*}}$ at $\p^{*}$ with \texttt{solver}}
		\vspace{0.5em}
		
		\STATE{Compute true defect $\{\mathbf{d}(t^{k}, \p)\}_{k=1}^{N_{t}}$ at $\p^{*}$ based on $\bX^{*}$}
		\vspace{0.5em}
		
		\STATE{Determine $\bV^{*}$ through an SVD of $\overline{\bX} := 
			\bX^{*} - \bV (\bV^{\tpose} \bX^{*})$}
		\vspace{0.5em}
		
		\STATE Update $\bV$ as $\bV := \texttt{orth}\big( \bV, \bV^{*}(: \,,\,1:r_{c})\big)$\\ with $\texttt{orth}\big(\cdot\big)$ denoting an orthogonalization process which can be implemented using the modified Gram-Schmidt process, or QR algorithm
		\vspace{0.5em}
		
		\STATE{Obtain the C-ROM (e.g., \cref{eq:rom_timedisc_corr}) corresponding to $\texttt{solver}_{\textnormal{imp}}$ by Galerkin proj. (+ hyperreduction)}
		\vspace{0.5em}
		
		\STATE{Solve the C-ROM to obtain the residual and compute the error estimator $\overline{\Delta}_{a}(\p)$ \cref{eq:data_enhanced_estimator_a}\\ or $\overline{\Delta}_{b}(\p)$ \cref{eq:data_enhanced_estimator}}
		\vspace{0.5em}
		
		\STATE{$\p^{*} := \arg \max \limits_{\p \in \Xi} \overline{\Delta}_{z}(\p)$, $z = a$ or $z = b$}
		\vspace{0.2em}
		
		\STATE{Set $\epsilon = \overline{\Delta}_{z}(\p^{*})$, $z = a$ or $z = b$}
		\vspace{0.5em}
		
		\ENDWHILE
	\end{algorithmic}
\end{algorithm}
\subsubsection{Computational cost}
We now analyse the additional computational cost incurred by the proposed algorithm \Cref{alg:podgreedy_ode} when compared to the standard POD-Greedy method in \Cref{alg:podgreedy}.

To simplify things, we define the cost of solving the nonlinear FOM \cref{eq:fom} to be  $\mathcal{C}_{\text{FOM}} := N_{\text{newton}} N_{L} N_{t}$; the cost of the linear solve at each Newton iteration, viz. $N_{L}$, could range between $N^{2} - N^{3}$ depending on the particular method being implemented in the solver. The cost of obtaining the defect vector at each time step (see \cref{eq:defect_im1}) is denoted by $\mathcal{C}_{\mathbf{d}}$. It depends on the user-defined time-stepping scheme. The main contributions within $\mathcal{C}_{\mathbf{d}}$  are matrix-vector products and evaluations of the nonlinearity. Thus, it evaluates to $\mathcal{C}_{\mathbf{d}} :=  ( \mathcal{O}(N^{2}) + \mathcal{O}(N))$ in the worst-case . However, the matrix vector multiplications typically involve sparse matrices and can be done cheaply. 

We denote by $\mathcal{C}_{\text{SVD}} := \mathcal{O}(\text{min}(p^{2}  q, p  q^{2}))$ the cost of a SVD for a matrix of dimension $\R^{p \times q}$. We let $\mathcal{C}_{\text{RBF}} := \mathcal{O}(d_{s}^{3})$ be the cost of obtaining one RBF interpolant. The cost of training a neural network is difficult to estimate owing to its architecture and the use of specialized hardware. Therefore, for simplicity, we define it as $\mathcal{C}_{\text{NN}} := \sigma \cdot (\mathcal{C}_{\text{FP}} + \mathcal{C}_{\text{BP}})$ with $\mathcal{C}_{\text{FP}}$ denoting the cost of one forward pass, $\mathcal{C}_{\text{BP}}$ denoting the cost of one backward pass and $\sigma$ being a constant that depends on the number of layers, epochs, the batch size and other hyperparameters.

The overall factors contributing to the proposed algorithm are listed below:

\begin{itemize}
	\item Step 2 in \Cref{alg:podgreedy_ode} requires the FOM solution at $d_{s}$ parameters. This has cost that scales as $d_{s} \cdot \mathcal{C}_{\text{FOM}}$.
	\item In Step 2, evaluating the defect vector at each time step and for all $d_{s}$ parameter samples incurs cost that scales as $(d_{s} N_{t}) \cdot \mathcal{C}_{\mathbf{d}}$.
	\item The cost of the two-stage SVD in Step 3 to obtain $\bV_{d}$ is $(d_{s}+1) \cdot \mathcal{C}_{\text{SVD}}$.
	\item Based on the method used for approximating the map $(t, \p) \mapsto \widehat{\bd}(t, \p)$, the costs differ:
		\begin{itemize}
			\item[-] For the RBF-based approach, the RBF coefficient matrix can be factorized once and reused for subsequent solves.  The cost is thus $\mathcal{C}_{\text{RBF}} + (n_{d} N_{t} - 1) \cdot 	\mathcal{O}(d_{s}^{2})$ 
			\item[-] For the NN-based approach, the cost is $\mathcal{C}_{\text{NN}}$
		\end{itemize}
	\item In the inference stage where the defect vector is approximated for all $\p \in \Xi$, 
		\begin{itemize}
			\item[-] the RBF-based approach has a cost $(N_{p} N_{t}) \cdot d_{s}$
			\item[-] the NN-based approach incurs a cost scaling as $(N_{p} N_{t}) \cdot \mathcal{C}_{\text{FP}}$
		\end{itemize}
	The cost of the matrix tensor product to obtain $\widetilde{\mathfrak{D}}$ is $d_{s} N_{t} \cdot n_{d} N $ for both RBF and NN-based approaches
	\item To update the approximate defect vector in Step 6, the cost of evaluating the true defect vector at the current greedy parameter scales as $N_{t} \cdot \mathcal{C}_{\mathbf{d}}$
\end{itemize}

The major cost in approximating the defect will be the cost of solving the FOM at $d_{s}$ parameter samples. From our experience, the RBF-based approach performs better than the NN-based approach. Comparisons of corresponding run times are detailed in \Cref{sec:numerics}.
\section{Numerical results}
\label{sec:numerics}
To demonstrate the validity of the proposed data-enhanced output error estimation approach, we test it on three numerical examples. These are:
\begin{enumerate}
	\item the viscous Burgers' equation with one parameter
	\item the FitzHugh-Nagumo equations with two parameters and
	\item the batch chromatography equations with one parameter.
\end{enumerate}
All reported numerical results were performed on a desktop computer running Ubuntu 20.04, installed with a $12$-th generation \INTEL \textsf{i5} processor, $32$GB of RAM and a NVIDIA RTX A4000 GPU with $16$GB of memory. The simulations were carried out in the Spyder IDE with Python 3.9.12 (\textsf{miniconda}). Where required, \MATLAB $\,$ v2019b was used to run some simulations. The radial basis interpolation was performed using the RBF package \cite{hinesrbf}.

For the two greedy algorithms (\Cref{alg:podgreedy,alg:podgreedy_ode}), we plot the maximum estimated errors computed using $\overline{\Delta}_{b}(\p)$ over the training set at every iteration. We define this as:
\[
	\varepsilon_{\text{max}} := \max \limits_{\p \in \Xi} \overline{\Delta}_{b}(\p)
\]
where $\overline{\Delta}_{b}(\p)$ is defined in \cref{eq:data_enhanced_estimator}. Additionally, to illustrate the performances of the ROM over the test set $\Xi_{\text{test}}$, we plot the mean estimated error $\overline{\Delta}_{b}(\p)$ for every parameter $\p \in \Xi_{\text{test}}$.

\subsection{Code availability}
The companion Python code to reproduce the numerical results is available at\\ \url{https://doi.org/10.5281/zenodo.8169490}.

\subsection{Burgers' equation}
\label{sec:numerics_burgers}
\paragraph{Model description}
The viscous Burgers' equation defined in the 1-D domain $\Omega := [0, 1]$ is given by
\begin{align}
\label{eq:burgers_pde}
\frac{\partial v}{\partial t} + v \frac{\partial v}{\partial z} &= \mu \frac{\partial^{2} v}{\partial^{2} z},\\ \nonumber
v(z, 0) &= \sin(2 \pi z), \\ \nonumber
v(0, t) &= v(1, t) = 0
\end{align}
with $v := v(z, t) \in \R^{}$ denoting the state variable and $z \in \Omega$ is the spatial variable and the time variable $t \in [0, 2]$. We spatially discretize \cref{eq:burgers_pde} with the finite difference method. The mesh size is $\Delta z = 0.001$, which results in a discretized FOM of dimension $N = 1000$. As the variable parameter, we consider the viscosity $\mu \in \mathcal{P} := [0.005, 1]$. We sample $100$ logarithmically-spaced samples from $\mathcal{P}$ and divide the samples randomly into a training set $\Xi$ and a testing set $\Xi_{\text{test}}$ in the ratio $80:20$. To solve the ODE, the \texttt{solver} in \Cref{alg:podgreedy_ode} is the \texttt{scipy} library \texttt{odeint}. To have a uniform comparison, we compute solutions of the FOM on a uniformly spaced time step of $\delta t = 0.01$. The output variable of interest is the value of the state at the node just before the right boundary.

\begin{figure}[t!]
	\centering
	\subfloat[\label{fig:burgers_imex1_none_convergence}]{\includegraphics[scale=0.7]{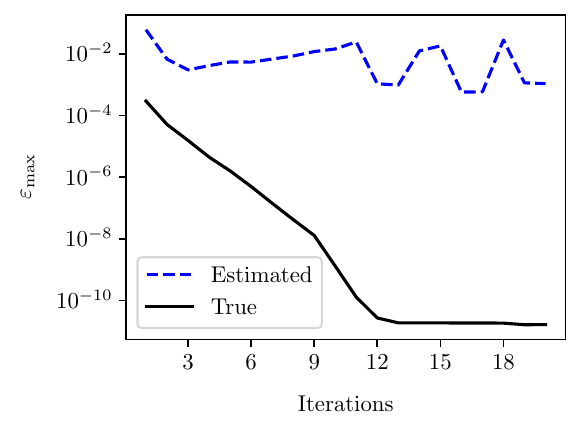}}
	\hspace{2em}
	\subfloat[\label{fig:burgers_imex1_true_convergence}]{\includegraphics[scale=0.7]{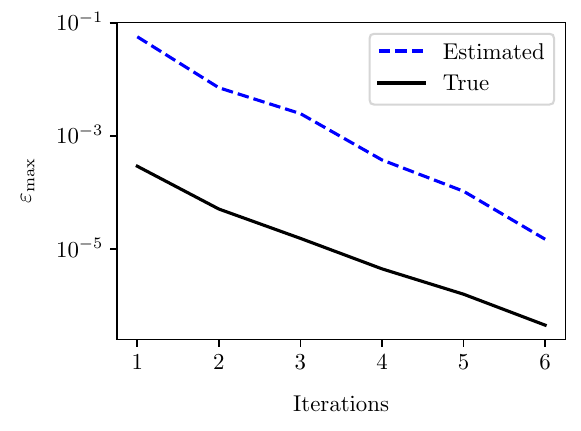}}
	\caption{Burgers' equation, \Cref{alg:podgreedy}: (a) error (estimator) decay when using a first-order IMEX method without any closure term; (b) error (estimator) decay when the true residual is known.}	
	\label{fig:burgers_imex1_none_true_convergence}%
\end{figure}

\begin{figure}[t!]
	\centering
	\subfloat[\label{fig:burgers_imex1_svdrbf_convergence}]{\includegraphics[scale=0.7]{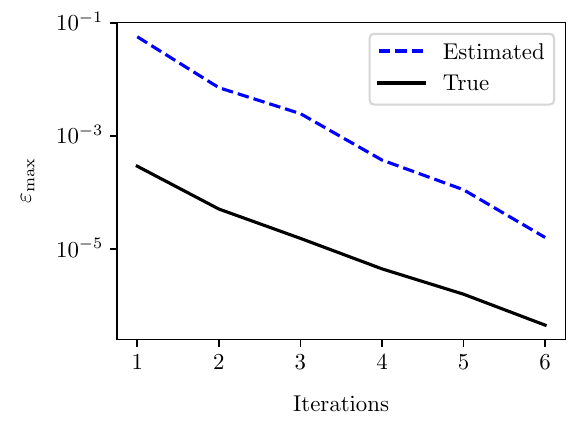}}
	\hspace{1em}
	\subfloat[\label{fig:burgers_imex1_svdrbf_testerror}]{\includegraphics[scale=0.7]{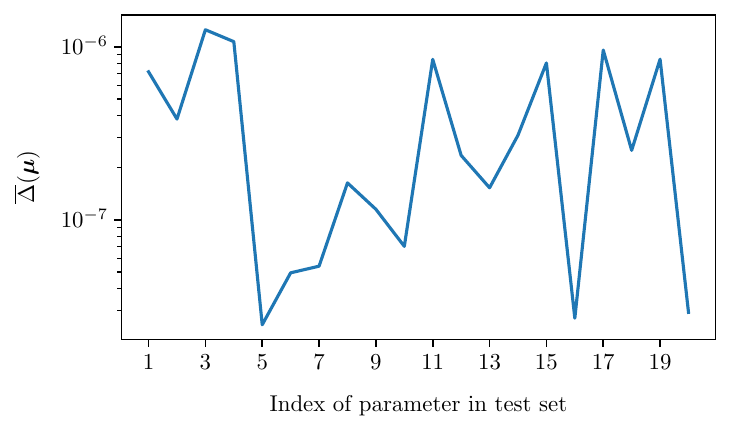}}
	\caption{Burgers' equation, \Cref{alg:podgreedy_ode}: (a) error (estimator) decay when using the SVD+RBF method to approximate the closure term and Step 6 is not included; (b) performance of the ROM over a test set.}
	\label{fig:burgers_imex1_svdrbf_convergence_testerror}%
\end{figure}

\begin{figure}[t!]
	\centering
	\subfloat[\label{fig:burgers_imex1_svdrbf_svd_t}]{\includegraphics[scale=0.7]{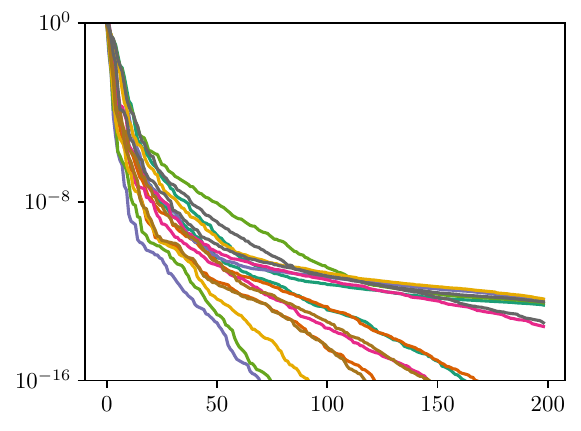}}
	\hspace{1em}
	\subfloat[\label{fig:burgers_imex1_svdrbf_svd_mu}]{\includegraphics[scale=0.7]{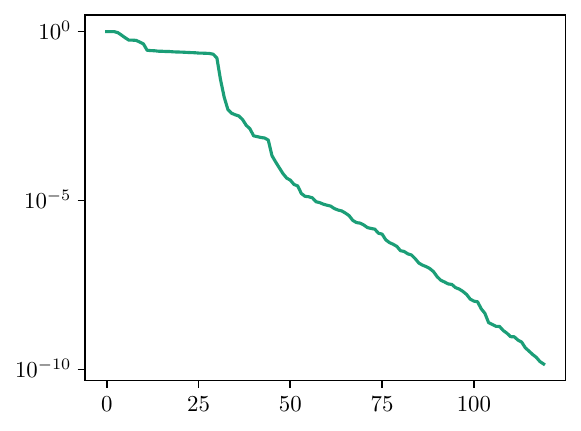}}
	\caption{Burgers' equation: (a) normalized singular values of the defect matrix $\mathbf{D}(\p)$ for $\p \in \Xi_{\text{defect}}$; (b) normalized singular values of $\mathbf{R}$.}
	\label{fig:burgers_imex1_svdrbf}%
\end{figure}

\paragraph{POD-Greedy with \Cref{alg:podgreedy}}
We first apply \Cref{alg:podgreedy} to the Burgers' equation with a tolerance $\texttt{tol} = 10^{-4}$. The $\texttt{solver}$ used is the $\texttt{odeint}$ library in \texttt{scipy}. Since the exact expression of the residual is unknown, we use a first-order IMEX method (IMEX1) to approximate the residual and estimate the output error. As the true residual is incorrectly approximated by the user-imposed time integration scheme, the resulting estimated error is severely overestimated. This results in the stagnation of the greedy algorithm as seen in \Cref{fig:burgers_imex1_none_convergence}. For comparison, we also plot in \Cref{fig:burgers_imex1_true_convergence}  the convergence of the greedy algorithm when the exact residual is known. To remedy this, we next apply the proposed \Cref{alg:podgreedy_ode}. 
\paragraph{POD-Greedy with \Cref{alg:podgreedy_ode} and SVD+RBF closure approximation}
First, we collect $d_{s} = 16$ uniformly-spaced parameter samples from the training set to construct $\Xi_{\text{defect}}$ and obtain the corresponding defect vector in Step 2. The SVD tolerances $\texttt{tol}_{\text{SVD}, t}, \texttt{tol}_{\text{SVD}, \p}$ are both set to $10^{-4}$ so that $n_{d} = 47$. The user-defined solver $\texttt{solver}_{\textnormal{imp}}$ is IMEX1.  The resulting convergence of the greedy algorithm using the RBF-based approximation of the defect vector is shown in \Cref{fig:burgers_imex1_svdrbf_convergence}. As shown, the maximum estimated error converges exponentially to the desired tolerance. The dimension of the ROM obtained is $n=7$. To demonstrate the performance of the ROM, we show in \Cref{fig:burgers_imex1_svdrbf_testerror} the mean estimated errors for the parameters in the test set $\Xi_{\text{test}}$. The obtained errors are smaller than the desired tolerance $10^{-4}$, showing the reliability of our error estimation approach. Furthermore, \Cref{fig:burgers_imex1_svdrbf_svd_t} shows the singular value decays of the defect vector trajectories at the $16$ parameter samples in $\Xi_{\text{defect}}$. \Cref{fig:burgers_imex1_svdrbf_svd_mu} shows the singular value decay obtained from the SVD of $\mathbf{R}$ (\Cref{subsec:SVD_defect_reduction}). The singular value decay in \Cref{fig:burgers_imex1_svdrbf_svd_mu} indicates the fast Kolmogorov $n$-width decay of the defect manifold and hence good performance of \Cref{alg:podgreedy_ode}. In terms of runtime, this approach requires $17$ seconds for generating the training data for learning the defect vector and obtaining the RBF interpolants.

\begin{figure}[b!]
	\centering
	\subfloat[\label{fig:burgers_imex1_svdnn_convergence_nodefectupdate}]{\includegraphics[scale=0.7]{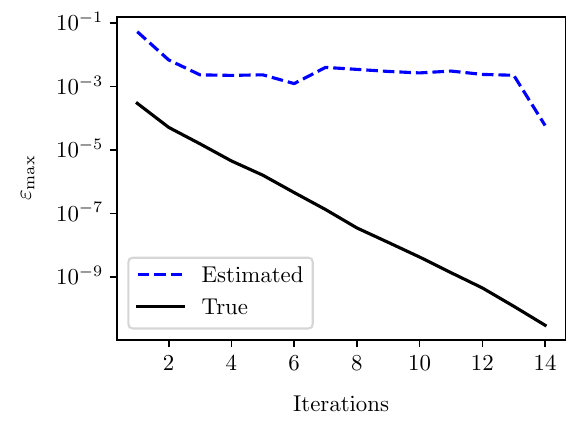}}
	\hspace{2em}
	\subfloat[\label{fig:burgers_imex1_svdnn_convergence_withdefectupdate}]{\includegraphics[scale=0.7]{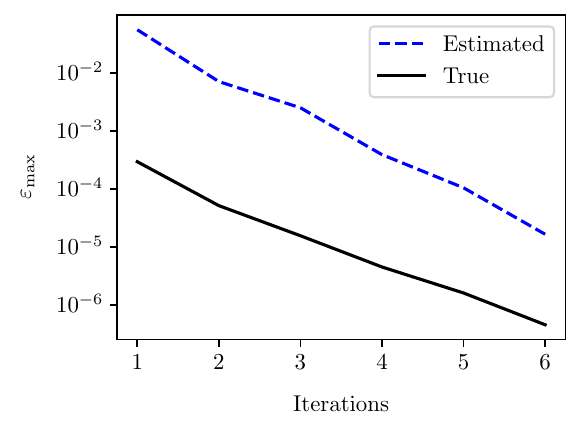}}
	\caption{Burgers' equation, \Cref{alg:podgreedy_ode}: (a) error (estimator) decay when using the SVD+FNN method without closure updates; (b) error (estimator) decay when using the SVD+FNN method with closure updates.}
	\label{fig:burgers_imex1_svdnn_convergence}%
\end{figure}

\paragraph{POD-Greedy with \Cref{alg:podgreedy_ode} and SVD+FNN closure approximation}
Next, we repeat \Cref{alg:podgreedy_ode} but now with the neural network-based approximation of the defect. The feed-forward neural network has 3 hidden layers with $16, 64, 64$ neurons, respectively. The activation function are the \textsf{SiLU} function for the first three layers and \textsf{Tanh} for the last layer. We have normalized the input data to be between $[0, 1]$ and the output data is between $[-1, 1]$. The neural network is implemented in PyTorch. It is trained using the \textsf{Adam} optimizer for $2000$ epochs, with the learning rate being set as $0.005$. Initially, we do not implement Step 6 and do not update $\widetilde{\bd}_{\text{NN}}(t, \p^{*})$ with $\mathbf{d}(t, \p^{*})$ at $\p^{*}$ (selected from the previous iteration) when computing the error estimator $\overline{\Delta}_{b}(\p)$ in Step 11 at the current iteration. Using the same $\Xi_{\text{defect}}$ with $|\Xi_{\text{defect}}| = 16$ (as done for the SVD+RBF approach), we did not obtain convergence of the greedy algorithm. Therefore, we use a $\Xi_{\text{defect}}$ with $24$ samples. The SVD tolerances for this case are both $0.1$ resulting in $n_{d} = 6$. Using enriched training data results in the successful convergence of \Cref{alg:podgreedy_ode} as seen in \Cref{fig:burgers_imex1_svdnn_convergence_nodefectupdate}. However, it requires up to $14$ iterations for this convergence. Evidently, the NN-based approach does not yield a satisfactory performance even with more training data. Then, we implement Step 6 and update the defect vector approximation at Step 11, which results in a  significant improvement in performance. This is shown in \Cref{fig:burgers_imex1_svdnn_convergence_withdefectupdate}.

\begin{figure}[t!]
	\centering
	\includegraphics[scale=0.55]{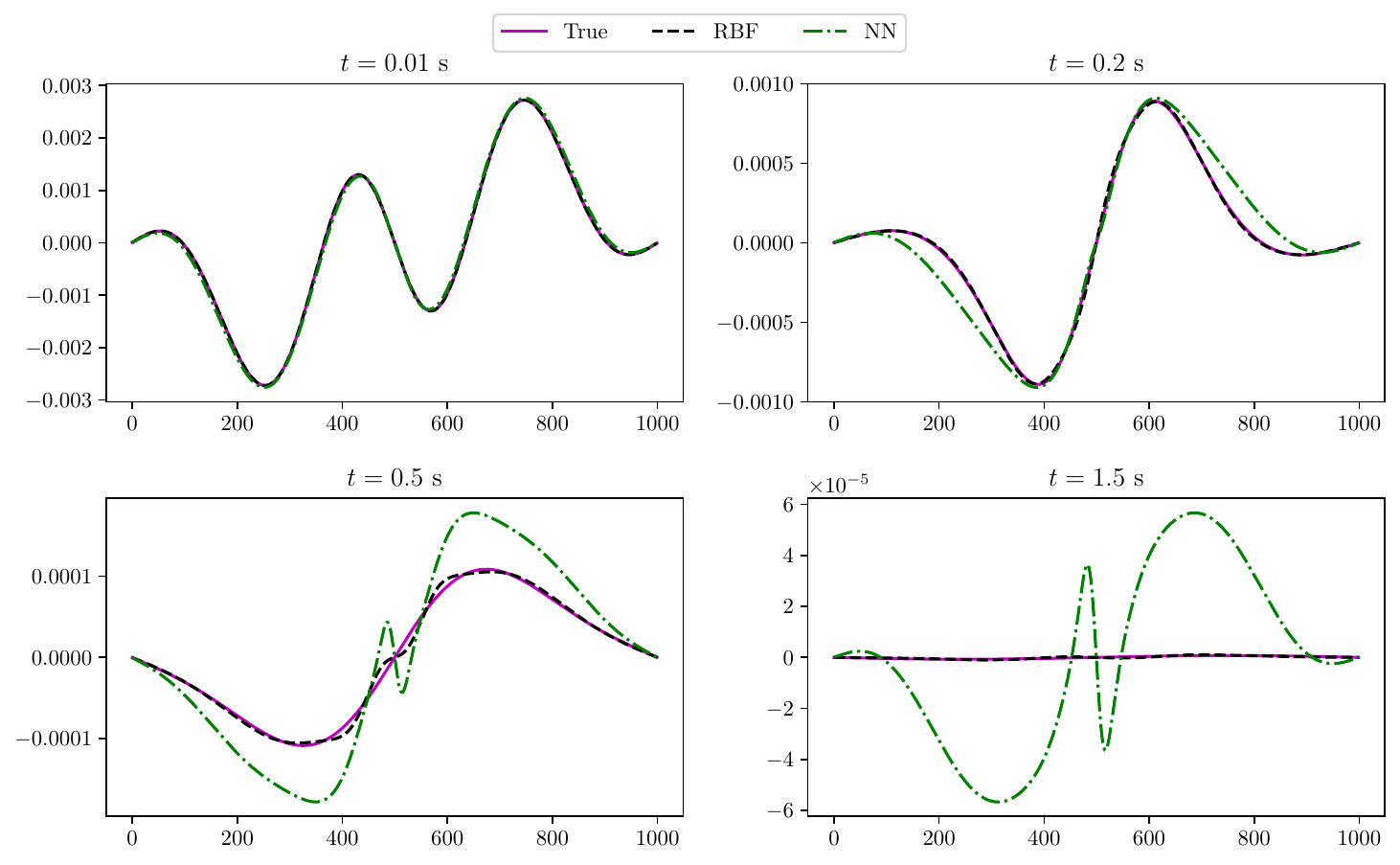}
	\caption{Burgers' equation: approximation of the true defect by the RBF-based and NN-based approaches at $t \in \{0.01, 0.2. 0.5, 1.5\}$ s for $\mu = 0.123$.}
	\label{fig:burgers_imex1_defectcompare}%
\end{figure}

To explain the poorer performance of the SVD+FNN approach, we plot in \Cref{fig:burgers_imex1_defectcompare} the approximation of the defect vector at $\mu = 0.123$ and $4$ different time instances, viz., $t \in \{0.01, 0.2, 0.5, 1.5\}$~s. We notice that the approximation from the FNN, while qualitatively capturing the true defect, fails to produce a very close match to the true value. This is especially the case for latter time instances, as the magnitude of the defect vector gets smaller. However, the RBF-based approach results in a significantly better approximation. This might be because of the fact that each entry of the reduced defect vector is separately learned by an individually-trained RBF interpolation, while the whole reduced defect vector is learned by a single and uniformly-trained FNN. However, if we use $n_{d}$ FNNs to learn the $n_{d}$ entries of the defect vector separately, the training will become much more expensive when $n_{d}$ is not very small.

\subsection{FitzHugh-Nagumo equations}
\label{sec:numerics_fhn}
\paragraph{Model description}
The FitzHugh-Nagumo system models the response of an excitable neuron or cell under an external stimulus. It finds applications in a variety of fields such as cardiac electrophysiology and brain modeling. The nonlinear coupled system of two partial differential equations defined in the domain $\Omega := [0, L]$ is given below:
\begin{subequations}
	\label{eq:fhn_pde}
	\begin{align}		
	\epsilon \frac{\partial v_{1}(z,t)}{\partial t} &= \epsilon^{2} \frac{\partial^{2} v_{1}(z,t)}{\partial z^{2}} + f(v_{1}(z,t)) - v_{2}(z,t) + c, \label{eq:fhn_a}\\
	\frac{\partial v_{2}(z,t)}{\partial t} &= b\, v_{1}(z,t) - \gamma v_{2}(z,t) + c, \label{eq:fhn_b}
	\end{align}		
\end{subequations}
with boundary conditions
\begin{align}
\label{eq:fhn_pde_bnd_cond}
\frac{\partial }{\partial z} v_{1}(0, t) &= -I_{\text{ext}}(t), \hspace{1em} \frac{\partial }{\partial z} v_{1}(L, t) = 0,
\end{align}
and initial conditions
\begin{align}
\label{eq:fhn_pde_init_cond}
v_{1}(z, 0) = 0.001, \hspace{1em} v_{2}(z, 0) = 0.001.
\end{align}
In the above equations, $v_{1}(z, t)$ and $v_{2}(z, t)$ represent the electric potential and the recovery rate of the potential, respectively.
The spatial variable is denoted by $z \in \Omega$ and the time $t \in [0, 5]$. The nonlinear term is represented by $f(v_{1}(z, t)) := v_{1} (v_{1} - 0.1) (1 - v_{1})$. The external stimulus is $I_{\text{ext}}(t) = 50000 t^{3} e^{-15t} $. The system has four free parameters $\epsilon, c, b, \gamma$. We fix $b=0.5$ and $\gamma=2$ while the two free parameters are $\p = [\epsilon, c] \in \mathcal{P} := [0.01, 0.04] \times [0.025, 0.075]$. A finite difference scheme is employed to spatially discretize \cref{eq:fhn_a,eq:fhn_b} with $512$ nodes used for each variable leading to a FOM of dimension $N = 1024$. We sample $100$ parameters uniformly from the domain $\mathcal{P}$ and randomly divide them into the training set $\Xi$ and the test set $\Xi_{\text{test}}$ in the ratio $70 : 30$. To solve the ODE, we use \texttt{ode15s} from \MATLAB. The time discretization is done on a uniform grid with $\delta t = 0.01$. The output variables of interest are the values of the two state variables at the node next to the leftmost  boundary.

We apply \Cref{alg:podgreedy,alg:podgreedy_ode} to the FitzHugh-Nagumo system. This is a particularly challenging example for both algorithms as the system exhibits a slow singular value decay. Of particular interest is the approximation of the limit cycle behaviour of the system for certain combinations of the two free parameters, $(\epsilon, c)$. The RBM tolerance is set as $\texttt{tol} = 10^{-3}$.

\begin{figure}[t!]
	\centering
	\includegraphics[scale=0.7]{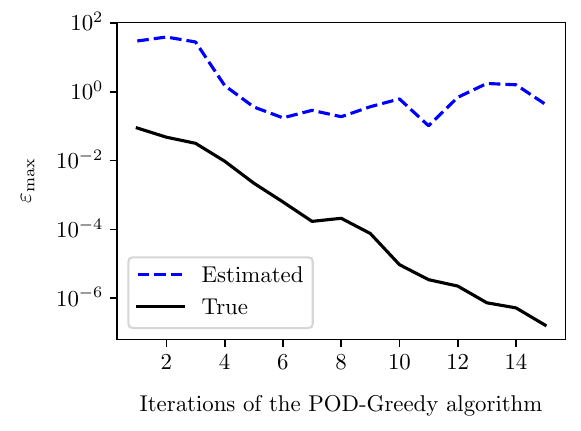}
	\caption{FitzHugh-Nagumo equations, \Cref{alg:podgreedy}: error (estimator) decay when using a second-order IMEX method without any closure term.}
	\label{fig:fhn_imex2_none_convergence}
\end{figure}

\paragraph{POD-Greedy with \Cref{alg:podgreedy}} Applying \Cref{alg:podgreedy} to this example and using a second-order IMEX scheme (IMEX2) to compute the residual does not result in the convergence of the greedy algorithm, as see in \Cref{fig:fhn_imex2_none_convergence}. This stems from the fact that the actual numerical scheme \texttt{solver} used in \Cref{alg:podgreedy} is the \texttt{ode15s} solver from \MATLAB, therefore, the residual we compute using IMEX2 scheme is incorrect.
\begin{figure}[t!]
	\centering
	\subfloat[\label{fig:fhn_imex2_svdrbf_svd_t}]{\includegraphics[scale=0.7]{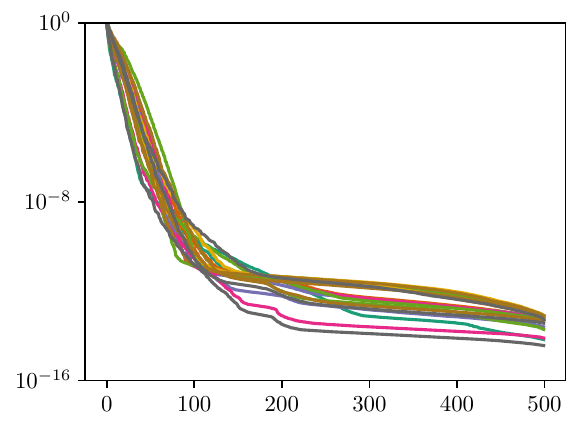}}
	\hspace{1em}
	\subfloat[\label{fig:fhn_imex2_svdrbf_svd_mu}]{\includegraphics[scale=0.7]{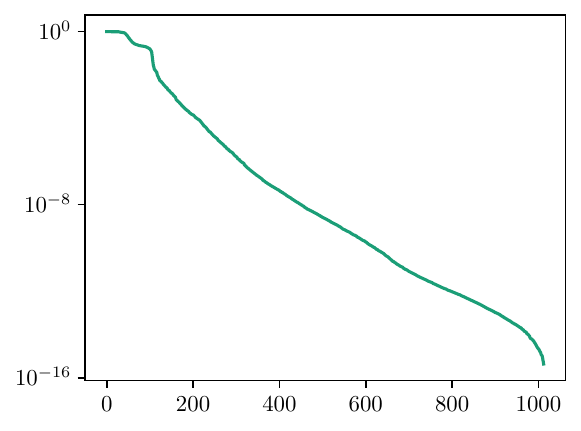}}
	\caption{FitzHugh-Nagumo equations: (a) Stage 1 - normalized singular values of the defect matrix $\mathbf{D}(\p)$ for $\p \in \Xi_{\text{defect}}$; (b)  Stage 2 - normalized singular values of $\mathbf{R}$.}
	\label{fig:fhn_imex2_svdrbf_svd_t_mu}%
\end{figure}
\begin{figure}[t!]
	\centering
	\subfloat[\label{fig:fhn_imex2_svdrbf_convergence_updatedefect}]{\includegraphics[scale=0.7]{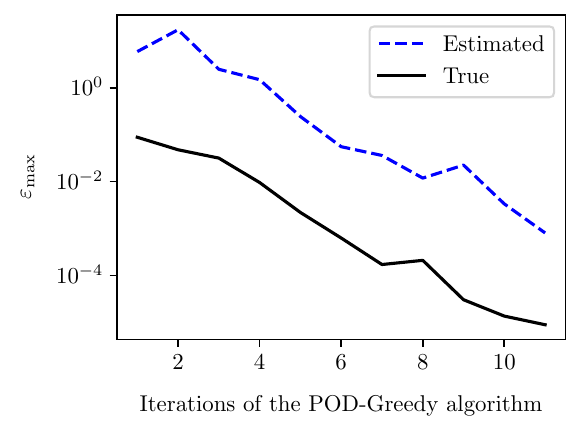}}
	\hspace{1em}
	\subfloat[\label{fig:fhn_imex2_svdrbf_convergence_noupdatedefect}]{\includegraphics[scale=0.7]{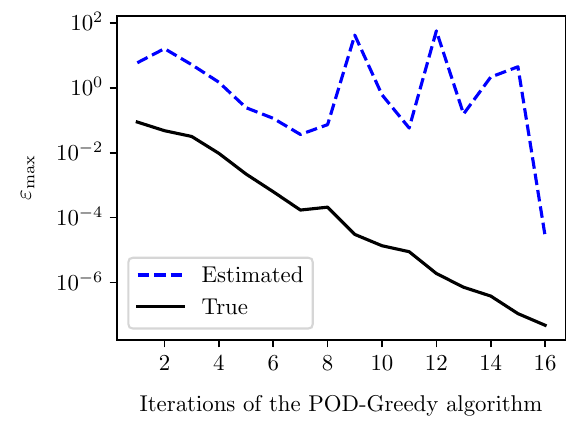}}
	\caption{FitzHugh-Nagumo equations, \Cref{alg:podgreedy_ode}: (a) error (estimator) decay when using the SVD+RBF method with Step 6 included; (b) error (estimator) decay when using the SVD+RBF method without Step 6.}
	\label{fig:fhn_imex1_svdrbf_convergence}%
\end{figure}

\paragraph{POD-Greedy with \Cref{alg:podgreedy_ode} and SVD+RBF closure approximation}
We apply \Cref{alg:podgreedy_ode} to this example using $\widetilde{\bd}_{\text{RBF}}(t, \p)$ in Step 3 and with an RBM tolerance $\texttt{tol} = 10^{-3}$. The user-imposed time integration scheme $\texttt{solver}_{\text{imp}}$ is IMEX2. Due to the challenging nature of the problem, $d_{s} = 21$ uniformly-spaced samples are chosen from the training set to obtain $\Xi_{\text{defect}}$. The tolerance $\texttt{tol}_{\text{SVD}, t} = \texttt{tol}_{\text{SVD}, \p} = 10^{-6}$, resulting in $n_{d} = 311$. \Cref{fig:fhn_imex2_svdrbf_svd_t} plots the singular value decays of $\mathbf{D}(\p)$ at all $\p \in \Xi_{\text{defect}}$ while \Cref{fig:fhn_imex2_svdrbf_svd_mu} shows the decay of the singular values of $\mathbf{R}$ (see \Cref{subsec:SVD_defect_reduction}). Similar to  the case of the Burgers' equation, an exponential decay of the singular values is observed in \Cref{fig:fhn_imex2_svdrbf_svd_t}. However, the second SVD shown in \Cref{fig:fhn_imex2_svdrbf_svd_mu} has a relatively slower decay of the singular values. This shows that the solution manifold of the FitzHugh-Nagumo system with respect to the parameter variations is more difficult to be approximated by a low-dimensional linear subspace. We incur $25$ seconds to compute the training data for the defect trajectories and to obtain the SVD+RBF approximation of the defect vectors. The greedy algorithm takes $11$ iterations to reach the desired tolerance; \Cref{fig:fhn_imex2_svdrbf_convergence_updatedefect} shows the error convergence of \Cref{alg:podgreedy_ode}. The dimension of the ROM obtained is $n=33$.
Note that we have implemented Step 6 in \Cref{alg:podgreedy_ode} to update the RBF approximation $\widetilde{\mathbf{d}}_{\text{RBF}}(t, \p^{*})$ with $\mathbf{d}(t, \p^{*})$ when we compute the error estimator $\overline{\Delta}_{b}(\p)$ for all $\p \in \Xi$ in Step 11. Without doing so, the greedy algorithm converges nevertheless, but takes $16$ iterations (\Cref{fig:fhn_imex2_svdrbf_convergence_noupdatedefect}) and the ROM has a larger size $n=48$.
\begin{figure}[t!]
	\centering
	\includegraphics[scale=0.7]{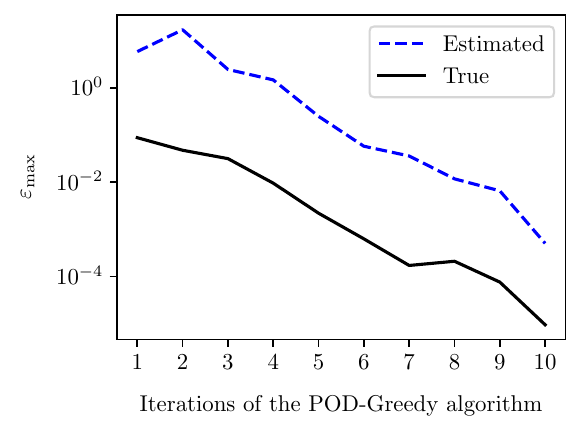}
	\caption{FitzHugh-Nagumo equations, \Cref{alg:podgreedy_ode}: error (estimator) decay when using the SVD+FNN method with Step 6 included.}
	\label{fig:fhn_imex2_svdnn_convergence}%
\end{figure}

\paragraph{POD-Greedy with \Cref{alg:podgreedy_ode} and SVD+FNN closure approximation}
Next, we use $\widetilde{\bd}_{\text{NN}}(t, \p)$ in Step 3 of \Cref{alg:podgreedy_ode}. The FNN is a $3$-layer network having, respectively, $64, 64, 32$ neurons in its hidden layers. The \textsf{SiLU} function is used for activation in all but the last layer. In the last layer, \textsf{Tanh} is the activation function. The training is carried out for $2000$ epochs using the \textsf{Adam} optimizer. The learning rate is $0.002$. No special tuning was done to calibrate the hyperparameters of the FNN. A detailed investigation on this is left for future work. We set $\texttt{tol}_{\text{SVD}, t} = \texttt{tol}_{\text{SVD}, \p} = 10^{-3}$, such that $n_{d} = 57$. The total time for computing the training data at all the samples in $\Xi_{\text{defect}}$ and for training the FNN is $106$ seconds, where training the FNN dominates the total runtime. \Cref{fig:fhn_imex2_svdnn_convergence} plots the convergence of the greedy algorithm. It takes $10$ iterations to converge. The resulting ROM has dimension $n=30$. 

\begin{figure}[t!]
	\centering
	\subfloat[\label{fig:fhn_imex2_testset_SVD+RBF_state}]{\includegraphics[scale=0.7]{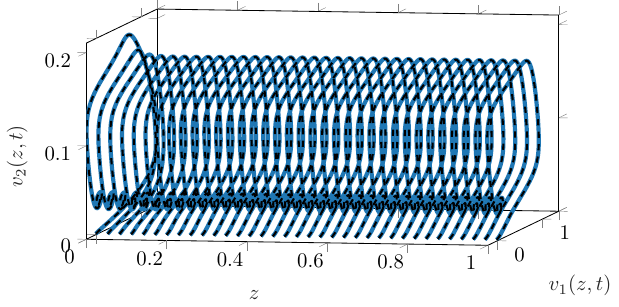}}
	\hspace{2em}
	\subfloat[\label{fig:fhn_imex2_testset_SVD+RBF_output}]{\includegraphics[scale=0.7]{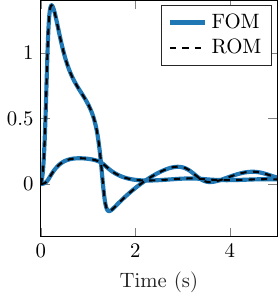}}
	\caption{FitzHugh-Nagumo equations, \Cref{alg:podgreedy_ode} with SVD+RBF: performance at the test parameter $\p = (0.0267, 0.0367)$ (a) Limit cycle behaviour; (b) output quantities.}
	\label{fig:fhn_imex2_testset_SVD+RBF}%
\end{figure}

\begin{figure}[t!]
	\centering
	\subfloat[\label{fig:fhn_imex2_testset_SVD+NN_state}]{\includegraphics[scale=0.7]{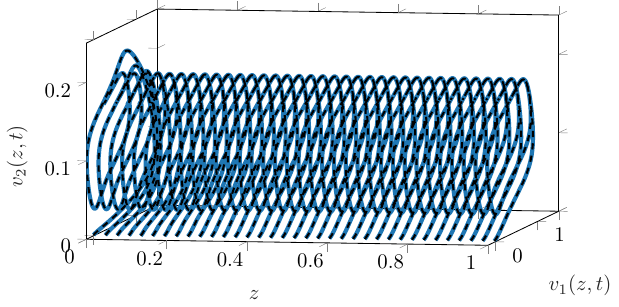}}
	\hspace{2em}
	\subfloat[\label{fig:fhn_imex2_testset_SVD+NN_output}]{\includegraphics[scale=0.7]{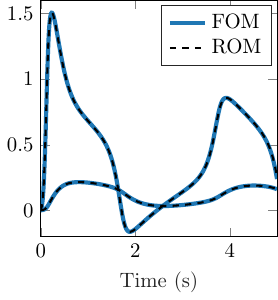}}
	\caption{FitzHugh-Nagumo equations, \Cref{alg:podgreedy_ode} with SVD+FNN: performance at the test parameter $\p = (0.04, 0.0472)$ (a) Limit cycle behaviour; (b) output quantities.}
	\label{fig:fhn_imex2_testset_SVD+NN}%
\end{figure}

The finally derived ROM is then simulated at two parameter samples $\p = (0.0267, 0.0367)$ and $\p = (0.04, 0.0472)$ taken from the test set. \Cref{fig:fhn_imex2_testset_SVD+RBF} shows the results of the ROM for the parameter $\p = (0.0267, 0.0367)$ obtained from \Cref{alg:podgreedy_ode} using the SVD+RBF approximation of the closure term. We see that the ROM is able to successfully capture both the state and the output dynamics of the FOM at the test parameter. At this parameter, the limit cycle behaviour is not very strong. At a different test parameter ($\p = (0.04, 0.0472)$) shown in \Cref{fig:fhn_imex2_testset_SVD+NN}, the ROM is able to successfully recover the stronger limit cycle behaviour as well. For this case, we show results using the SVD+FNN approach. But, we note that similar accuracy is also achieved with the SVD+RBF approach. The corresponding results are not shown to avoid repetition.

\subsection{Batch chromatography}
\label{sec:numerics_batchchrom}
\begin{figure}[t!]
	\centering
	\includegraphics[scale=0.25]{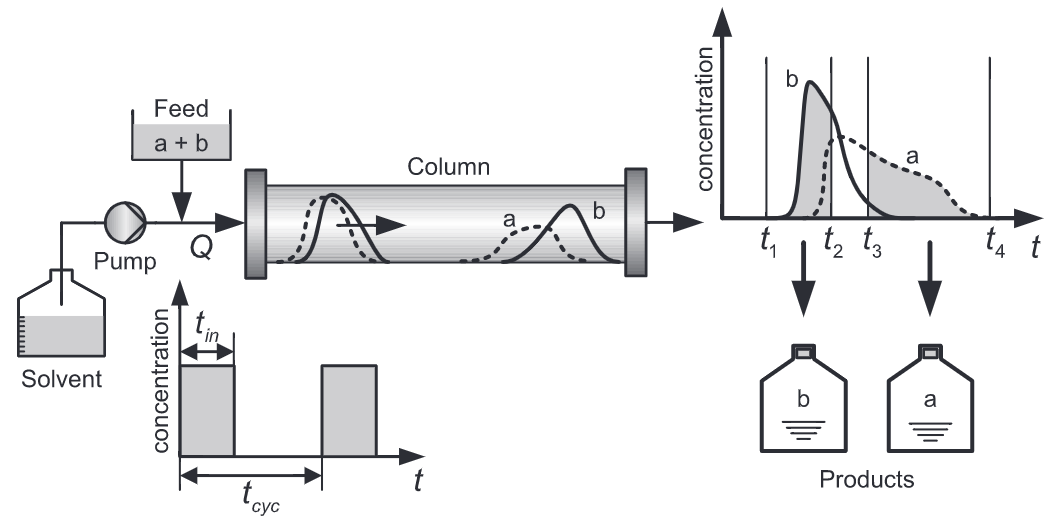}
	\caption{Schematic of the batch chromatography process.}
	\label{fig:batchchrom_schematic}
\end{figure}

The last example we consider is the model of the batch chromatography purification process. Being a coupled, nonlinear system of four PDEs, this example poses a considerable challenge for our proposed approach.
\paragraph{Model description}
Batch chromatography is an important purification process for separation of chemicals in food and pharmaceutical industries. We consider the governing equations for the batch chromatographic process for binary separation, i.e., the separation of two chemical components from a mixture. A schematic of the complete process is shown in \Cref{fig:batchchrom_schematic}.

\begin{figure}[t!]
	\centering
	\includegraphics[scale=0.7]{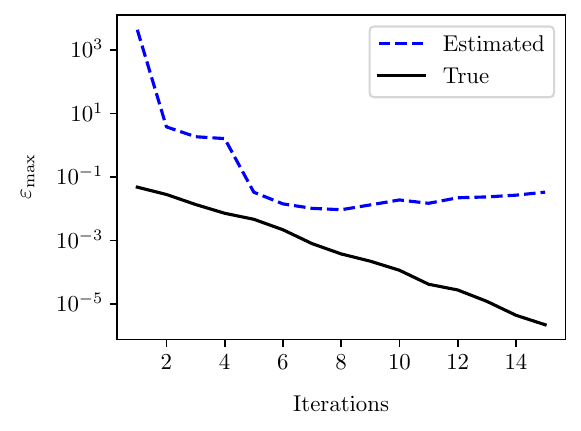}
	\caption{Batch chromatography, \Cref{alg:podgreedy}: error (estimator) decay when using a second-order IMEX method without any closure term.}
	\label{fig:batchchrom_imex2_none_convergence}
\end{figure}

The governing PDEs for the batch chromatography system are:
\begin{align}
\label{eq:batchchrom_pde}
\frac{\partial v_{1, z}}{\partial t} + \frac{1 - \epsilon}{\epsilon} \frac{\partial v_{2, z}}{\partial t} &= -\frac{\partial v_{1, z}}{\partial z} + \frac{1}{\text{Pe}} \frac{\partial^{2} v_{1,z}}{\partial z^{2}}, \\
\frac{\partial v_{2, z}}{\partial t} &= \dfrac{L}{Q/\epsilon A_{c}} \kappa_{z} \big(v_{2, z}^{\text{\scriptsize{Eq}}} - v_{2, z}\big)
\end{align}
where the state variables $v_{1,z}, v_{2,z}$ refer to the concentrations of the chemical component $z$ in the liquid and solid phase, respectively. Since we are interested in binary separation, $z \in \{a, b\}$. 	
The boundary conditions are:
\begin{align}
\label{eq:batchchrom_pde_bnd_cond}
\frac{\partial v_{1, z}(0, t)}{\partial z}  = \text{Pe} \big(v_{1,z}(0, t) - u(t)\big), \hspace{1em} \frac{\partial v_{1, z}(1,t)}{\partial z} = 0
\end{align}
and the initial conditions are:
\begin{align}
\label{eq:batchchrom_pde_init_cond}
v_{1,z}(z, 0) = 0, \hspace{1em} v_{2,z}(z, 0) = 0.
\end{align}
The quantity $v_{2, z}^{\text{\scriptsize{Eq}}}$ in \cref{eq:batchchrom_pde} is the source of nonlinearity and it denotes the adsorption equilibrium:
\begin{align*}
v_{2, z}^{\text{\scriptsize{Eq}}} = f_{z}(v_{1,a}, v_{2, b}) := \dfrac{H_{z1} v_{1, z}}{1 + K_{a1} v_{1,a}^{\text{f}} v_{1, a} + K_{b1} v_{2,b}^{\text{f}} v_{2, b}} +
\dfrac{H_{z2} v_{1, z}}{1 + K_{a2} v_{1,a}^{\text{f}} v_{1, a} + K_{b2} v_{2,b}^{\text{f}} v_{2, b}}.
\end{align*}
The discretization of the above PDE is performed using second-order finite volume method. Each of the four PDEs is spatially discretized into $800$ volume elements, resulting in a FOM of dimension $N=3200$. For full details regarding the batch chromatography PDE discretization and the terms involved, we refer to \cite{morChe23}. The output variables are the concentrations of the liquid phases ($v_{1,a}, v_{1,b}$) at the rightmost node. The batch chromatography model has two free parameters, $Q$ and $t_{\text{in}}$, which denote, respectively, the volumetric feed flow of the solvent injected into the column and the injection frequency of the solvent. We fix $t_{\text{in}} = 0.5$ while $Q \in [0.0667, 0.1667]$. We collect $60$ uniformly-spaced samples of $Q$ and divide them into the training set $\Xi$ and the test set $\Xi_{\text{test}}$ in the ratio $80 : 20$. The \texttt{solver} used for both algorithms is \texttt{ode15s} from \MATLAB while the user-imposed time integration method $\texttt{solver}_{\text{imp}}$ in \Cref{alg:podgreedy_ode} is IMEX2. The time discretization divides the time range $t \in [0, 10]$ into a uniform grid with $\delta t = 0.005$.
\paragraph{POD-Greedy with \Cref{alg:podgreedy}} First, we show the results of applying \Cref{alg:podgreedy} to the batch chromatography equations. The RBM tolerance is set to be $\texttt{tol} = 10^{-3}$. Since the exact form of the residual expression of batch chromatography equations is unknown, (i.e., $\mathcal{R}[\cdot]$ corresponding to \texttt{ode15s} is not available) we impose IMEX2 to get an approximate residual operator $\widetilde{\mathcal{R}}[\cdot]$. However, $\widetilde{\mathcal{R}}[\cdot]$  corresponding to IMEX2 is different from $\mathcal{R}[\cdot]$ so that the estimated error is inaccurate, resulting in stagnation of the greedy algorithm, see \Cref{fig:batchchrom_imex2_none_convergence}. To address this situation, we next apply the proposed approach which involves adding a closure term.

\begin{figure}[t!]
	\centering
	\includegraphics[scale=0.7]{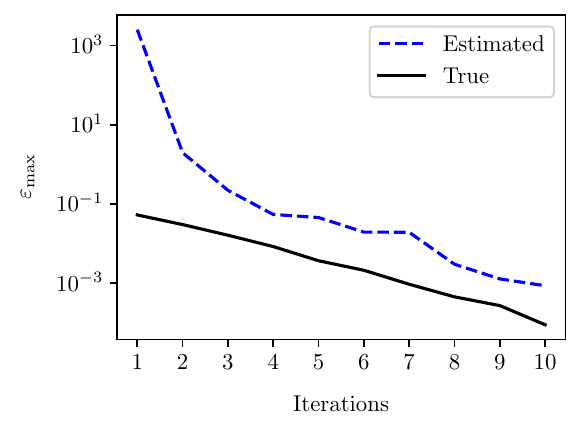}
	\caption{Batch chromatography, \Cref{alg:podgreedy_ode}: convergence when using the SVD+RBF method with Step 6 included to approximate the closure term.}
	\label{fig:batchchrom_imex2_svdrbf_convergence}%
\end{figure}

\begin{figure}[t!]
	\centering
	\subfloat[\label{fig:batchchrom_imex2_svdrbf_svd_t}]{\includegraphics[scale=0.7]{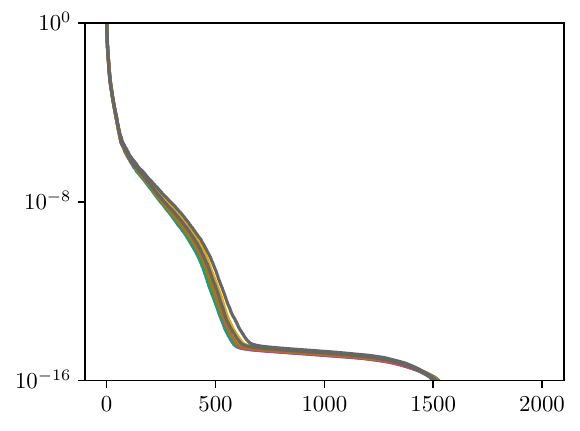}}
	\hspace{1em}
	\subfloat[\label{fig:batchchrom_imex2_svdrbf_svd_mu}]{\includegraphics[scale=0.7]{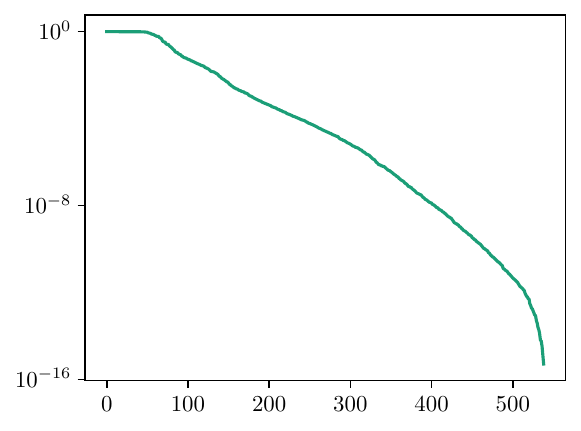}}
	\caption{Batch chromatography: (a) Stage 1 - normalized singular values of the defect matrix $\mathbf{D}(\p)$ for $\p \in \Xi_{\text{defect}}$; (b)  Stage 2 - normalized singular values of $\mathbf{R}$.}
	\label{fig:batchchrom_imex2_svdrbf_svd_t_mu}%
\end{figure}

\begin{figure}[t!]
	\centering
	\includegraphics[scale=0.8]{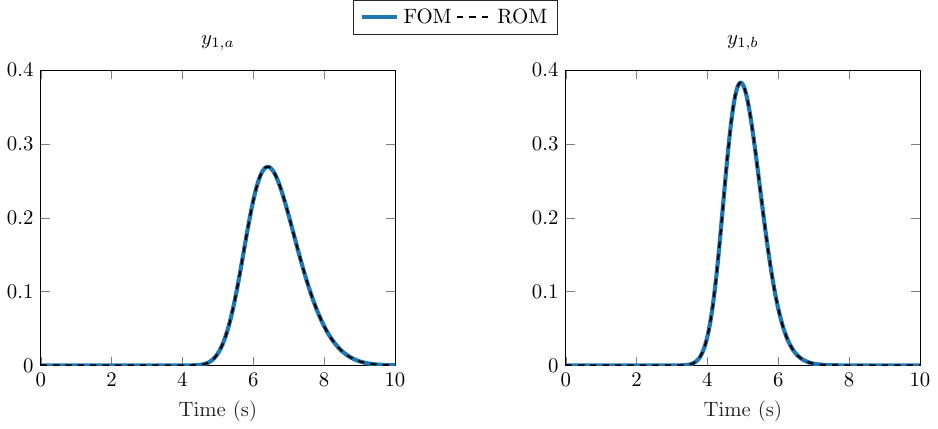}
	\caption{Batch chromatography, \Cref{alg:podgreedy_ode}: performance of the ROM at test parameter sample $Q=0.0803$ for the two output quantities $y_{1,a}, y_{1,b}$.}
	\label{fig:batchchrom_imex2_svdrbf_FOMROM_output}%
\end{figure}

\begin{figure}[t!]
	\centering
	\includegraphics[scale=0.8]{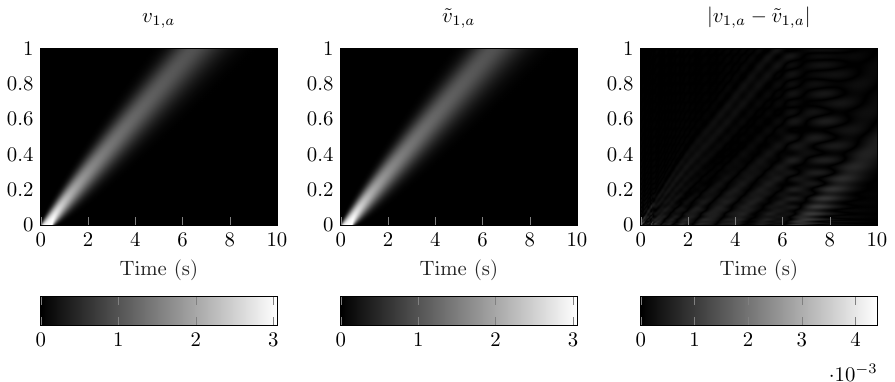}
	\caption{Batch chromatography, \Cref{alg:podgreedy_ode}: performance of the ROM at test parameter sample $Q=0.0803$ for the state $v_{1, a}$.}
	\label{fig:batchchrom_imex2_svdrbf_FOMROM_x1a}%
\end{figure}

\begin{figure}[t!]
	\centering
	\includegraphics[scale=0.8]{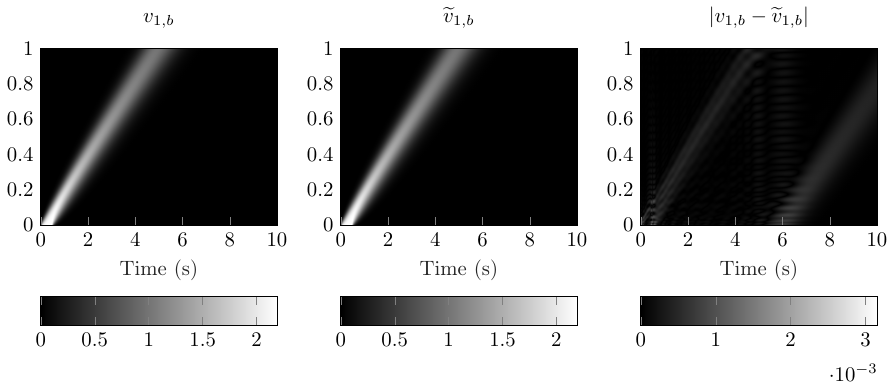}
	\caption{Batch chromatography, \Cref{alg:podgreedy_ode}: performance of the ROM at test parameter sample $Q=0.0803$ for the state $v_{1, b}$.}
	\label{fig:batchchrom_imex2_svdrbf_FOMROM_x1b}%
\end{figure}

\paragraph{POD-Greedy with \Cref{alg:podgreedy_ode} and SVD+RBF closure approximation}
As done for the last two examples, we first apply \Cref{alg:podgreedy_ode} to the current example and employ the RBF-based approach to approximate the closure term. The greedy algorithm successfully converges in $10$ iterations to the desired tolerance of $\texttt{tol} = 10^{-3}$. The dimension of the ROM is $n=134$. The convergence is shown in \Cref{fig:batchchrom_imex2_svdrbf_convergence}. To learn the closure term, $d_{s} = 10$ uniformly-spaced samples are chosen from the training set $\Xi$. The SVD tolerances are $\texttt{tol}_{\text{SVD}, t} = \texttt{tol}_{\text{SVD}, \p} = 10^{-5}$, leading to $n_{d} = 291$. The singular value decay from the SVD of each defect snapshot matrix $\mathbf{D}(\p)$ for all $\p \in \Xi_{\text{defect}}$ is presented in \Cref{fig:batchchrom_imex2_svdrbf_svd_t}. It can be noticed that the singular values decay even slower than those of the FitzHugh-Nagumo model (see \Cref{fig:fhn_imex2_svdrbf_svd_t}). In \Cref{fig:batchchrom_imex2_svdrbf_svd_mu}, the singular values of $\mathbf{R}$ also exhibit a much slower decay. This indicates that for the batch chromatography example, both the dynamics at a given parameter and the solution manifold with respect to the parameter variations are much more difficult to be captured by a low-dimensional linear space. As a result, this problem is likely to have a slow Kolmogorov $n$-width decay. The batch chromatography equations are in fact a system of first-order hyperbolic PDEs~\cite{LeeS19}. It is known that for problems exhibiting hyperbolic characteristics or convection-dominance, the Kolmogorov $n$-width decay is slow~\cite{GreU19}. 
The runtime for obtaining the training data for the SVD+RBF approach and for obtaining the RBF interpolant is $186$ seconds. To demonstrate the quality of the ROM resulting from \Cref{alg:podgreedy_ode}, we compare the output obtained from the FOM and ROM in \Cref{fig:batchchrom_imex2_svdrbf_FOMROM_output}. The results are shown for the parameter sample $Q = 0.0803$ taken from $\Xi_{\text{test}}$. We can see that the ROM is able to successfully recover the dynamics of both output quantities. For this sample, we obtained the mean error $\overline{\Delta}(Q) = 2.082 \cdot 10^{-5}$, which is below the desired tolerance. Additionally, to show the approximation of the state vector, we plot the space-time values and corresponding approximation errors of the liquid phase concentration $v_{1,a}, v_{1,b}$. The results are shown in \Cref{fig:batchchrom_imex2_svdrbf_FOMROM_x1a,fig:batchchrom_imex2_svdrbf_FOMROM_x1b}. It can be inferred that the ROM delivers a sharp approximation of both these state quantities for the entire duration of the simulation at an unseen parameter during training. Note that the state vector has error larger than the tolerance. The reason for this is that our error estimator aims to estimate the output error rather than the whole state error.

\begin{figure}[t!]
	\centering
	\includegraphics[scale=0.55]{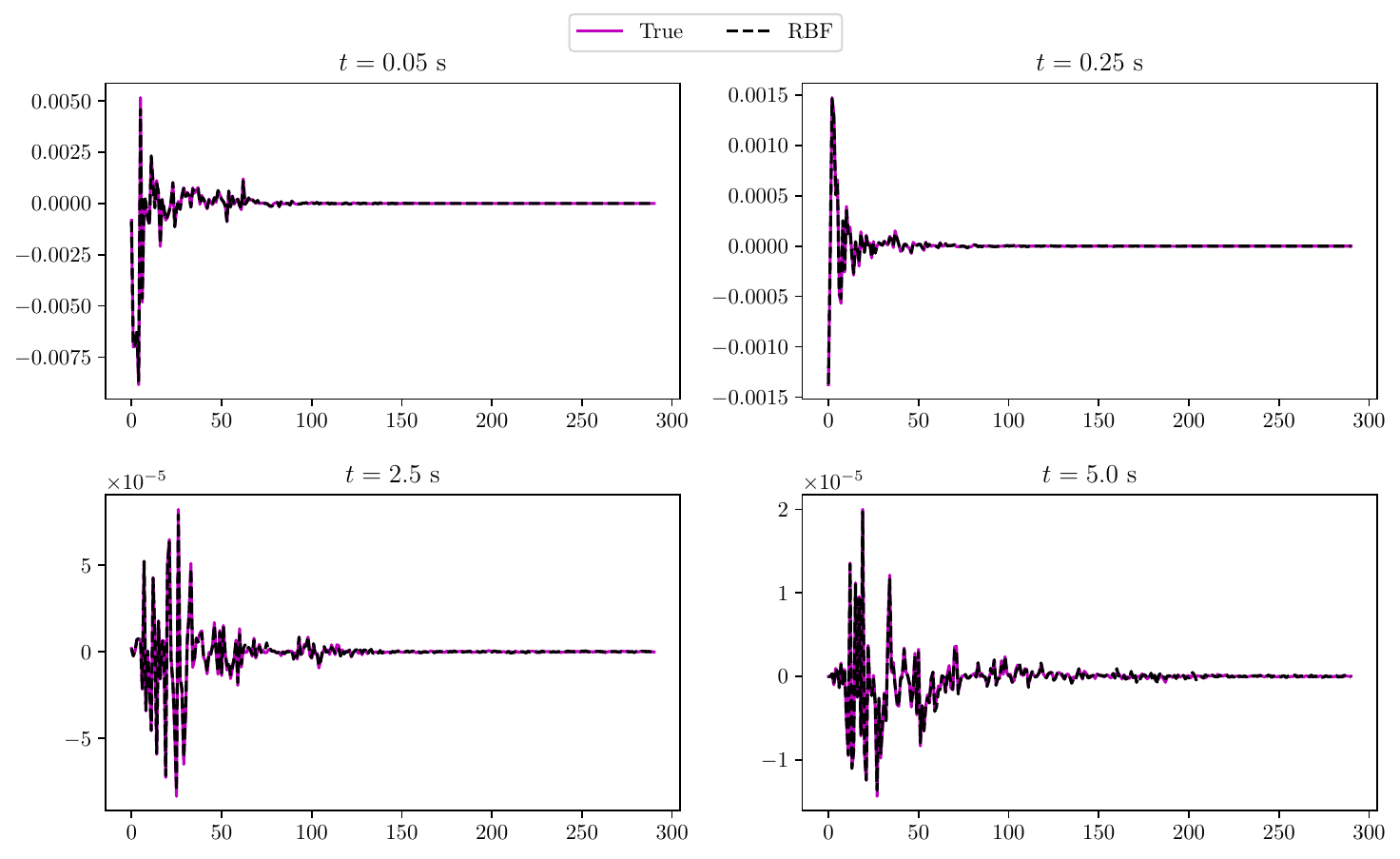}
		\caption{Batch chromatography: approximation of the reduced defect vector by the RBF-based approach at $t \in \{0.05, 0.25, 2.5, 5.0\}$ s at $Q = 0.0769$.}
	\label{fig:batchchrom_imex2_defectcompare}%
\end{figure}

The quality of the defect vector approximation using the SVD+FNN approach was not satisfactory for this example. This owes to the particularly non-smooth nature of the defect snapshots for different time instances and parameters. \Cref{fig:batchchrom_imex2_defectcompare} shows the true reduced defect vector $\widehat{\bd}(t, \p)$ and its approximation using RBF interpolation $\widehat{\bd}_{\text{RBF}}(t, \p)$ corresponding to $\p = 0.0769$ at the time instances $t \in \{0.05, 0.25, 2.5, 5.0\}$ s. While the RBF interpolants recover a good approximation, the FNN was not successful in capturing the entire complexity of the defect snapshots. The FNN-based approximation was very inaccurate and we do not show those results. The neural network struggles to capture the fast changing nature of the reduced defect vector in the reduced coordinate space and also its wide range of magnitudes ($\pm(10^{-3} - 10^{-8})$). The fact that we used separate RBF interpolants for each coordinate and time instance led to a much better approximation than the FNN. 
%
%
%
\section{Conclusion}
\label{sec:conclusion}
In this work, we introduced a data-enhanced \emph{a posteriori} output error estimator for model reduction of general parametric nonlinear dynamical systems. The proposed error estimator does not require any knowledge of the underlying time integration scheme used to integrate the given ODE. Applied to the reduced basis method, the new approach enables the direct use of ODE solver libraries, a feature that was not considered so far, to the best of our knowledge. While it demands a modest amount of extra training data, the proposed error estimator is efficient and can also be used in the online stage to certify the accuracy of the ROMs. Numerical experiments performed on three challenging examples demonstrate the benefits offered by the new approach. We observed that the RBF-based approach performed better, compared to a NN-based approach. An immediate extension of the proposed approach is to consider an adaptive sampling of the training set $\Xi$. New strategies to accurately approximate the closure term need to be considered for this. Another fruitful line of future work could involve integrating the proposed methodology in the Neural ODE framework~\cite{Cheetal18} to make it fully non-intrusive.

\section*{Acknowledgments}%
\addcontentsline{toc}{section}{Acknowledgments}
Part of this work was performed while the first author was pursuing his doctoral study and was supported by the International Max Planck Research School in Process Systems Engineering (IMPRS-ProEng).
\begin{appendices}
	\crefalias{section}{appendix}
	\section{Appendix A}
	\label{app:appendixA}
	\subsection*{Proof of \Cref{thm:resd_err_bound_state}}
	\begin{proof}
		Consider the FOM in \cref{eq:fom_timedisc} and its residual \cref{eq:resd_im1} based on the ROM \cref{eq:rom_timedisc}. Subtracting \cref{eq:resd_im1} from \cref{eq:fom_timedisc} and rearranging yields
		\begin{align}
		\bE_{\text{im}} \mathbf{e}^{k} &= \bA_{\text{im}} \mathbf{e}^{k-1} + \delta t \bigg(\bff(\bx_{\text{im}}^{k-1}) - \bff(\btx_{\text{im}}^{k-1})\bigg) + \br^{k}, \nonumber \\
		\mathbf{e}^{k} &= \bE_{\text{im}}^{-1} \bA_{\text{im}} \mathbf{e}^{k-1} + \delta t \bE_{\text{im}}^{-1} \bigg(\bff(\bx_{\text{im}}^{k-1}) - \bff(\btx_{\text{im}}^{k-1})\bigg) + \bE_{\text{im}} ^{-1} \br^{k}
		\end{align}
		where $\mathbf{e}^{k} = \| \bx_{\text{im}}^{k} - \btx_{\text{im}}^{k} \|$ is the error in the state vector at the $k$-th time step.
		
		Taking the norm on both sides and enforcing the Lipschitz condition results in
		\begin{align}
		\label{eq:state_error_recurrence}
		\nrm[2]{\mathbf{e}^{k}} &\leq \nrm[2]{\bE_{\text{im}}^{-1} \bA_{\text{im}}} \nrm[2]{\mathbf{e}^{k-1}} + \delta t L_{\bff} \nrm[2]{\bE_{\text{im}}^{-1}} \nrm[2]{\mathbf{e}^{k-1}}  + \nrm[2]{\bE_{\text{im}} ^{-1}} \nrm[2]{\br^{k}} \nonumber \\	
		&\leq \bigg( \nrm[2]{\bE_{\text{im}}^{-1} \bA_{\text{im}}} +
		\delta t L_{\bff} \nrm[2]{\bE_{\text{im}}^{-1}} \bigg)\nrm[2]{\mathbf{e}^{k-1}} +
		\nrm[2]{\bE_{\text{im}} ^{-1}}\nrm[2]{\br^{k}}.
		\end{align}
		For notational convenience, we define $\zeta := \nrm[2]{\bE_{\text{im}} ^{-1}}$ and $\xi := \big(\nrm[2]{\bE_{\text{im}}^{-1} \bA_{\text{im}}} +
		\delta t L_{\bff} \nrm[2]{\bE_{\text{im}}^{-1}}\big)$ and rewrite \cref{eq:state_error_recurrence} as
		\begin{align}
		\label{eq:state_error_recurrence_short}
		\nrm[2]{\mathbf{e}^{k}} &\leq \xi \nrm[2]{\mathbf{e}^{k-1}} +
		\zeta \nrm[2]{\br^{k}}.
		\end{align}
		
		At $k=0$, the initial error is 
		\begin{align}
		\label{eq:initial_proj_error}
		\nrm[2]{\mathbf{e}^{0}} = \nrm[2]{\bx_{\text{im}}^{0} - \bV \bV^{\tpose} \bx_{\text{im}}^{0}}.
		\end{align}
		
		Using \cref{eq:initial_proj_error}, the recursion in \cref{eq:state_error_recurrence_short} can be resolved for every $k$ to obtain the expression for the error bound in \cref{eq:std_errestm}.
	\end{proof}

%
	\section{Appendix B}
	\label{app:appendixB}
	\subsection*{Proof of \Cref{thm:err_estm}}
	\begin{proof}
		The output error for the modified output in \cref{eq:modified_output} is
		\begin{align} \label{eq:t_ymod}
		\by^{k} - \overline{\by}_{\text{im,c}}^{k} = \bC \big(\bx^{k} - \btx_{\text{im,c}}^{k} \big) + \btx_{\text{du}}^{\tpose} \br_{\text{im,c}}^{k}.
		\end{align}
		Multiplying by $\big( \bx^{k} - \btx_{\text{im,c}}^{k} \big)^{\tpose}$ on both sides of \cref{eq:fom_im1_dual} yields
		\begin{align*}
		\big( \bx^{k} - \btx_{\text{im,c}}^{k} \big)^{\tpose} \bE_{\text{im}}^{\tpose} \bxd &= - \big( \bx^{k} - \btx_{\text{im,c}}^{k} \big)^{\tpose} \bC^{\tpose},
		\end{align*}
		where we have made use of the fact that $\bE_{\text{du}} = \bE_{\text{im}}^{\tpose}$. Taking the transpose on both sides of the above equality leads to
		\begin{align}
		\label{eq:resd_output_reln}
		\bxd^{\tpose} \bE_{\text{im}} \big(\bx^{k} - \btx_{\text{im,c}}^{k} \big) &= -\bC \big(\bx^{k} - \btx_{\text{im,c}}^{k}\big).
		\end{align}
		We recall the auxiliary residual introduced in \cref{eq:resdbr_im1_corr} which can be written as
		\begin{align}
		\label{eq:resdbr_im1_corr_recall}
			\breve{	\br}_{\text{im,c}}^{k} &= \bA_{\text{im}} \bx_{\text{im,c}}^{k-1} + \delta t \big(\bff(\bx_{\text{im,c}}^{k-1}) + \bB \bu^{k} \big) + \bd^{k} - \bE_{\text{im}} \btx_{\text{im,c}}^{k},\nonumber\\
			&= \bE_{\text{im}} \big(\bx_{\text{im,c}}^{k} - \btx_{\text{im,c}}^{k}\big),\nonumber\\
			&=  \bE_{\text{im}} \big(\bx^{k} - \btx_{\text{im,c}}^{k}\big).
		\end{align}
		The third equality above follows from \cref{eq:defect_im1,eq:fom_timedisc_corr} that the corrected solution $\bx_{\text{im,c}}^{k}$ actually recovers the FOM solution $\bx^{k}$ given an accurate closure term $\bd^{k}$.
		
		We further use the expression in \cref{eq:resdbr_im1_corr_recall} to write \cref{eq:resd_output_reln} as
		\begin{align}
		\label{eq:t_err1}
		\bxd^{\tpose} \breve{\br}_{\text{im,c}}^{k} = -\bC \big(\bx^{k} - \btx_{\text{im,c}}^{k}\big).
		\end{align}
		Now, we substitute \cref{eq:t_err1} into \cref{eq:t_ymod} followed by addition and subtraction of the term $\btxdu^{\tpose} \breve{\br}_{\text{im,c}}^{k}$ to get
		\begin{equation}
		\label{eq:t_ymod2}
		\begin{aligned}
		\by^{k} - \overline{\by}_{\text{im,c}}^{k} &= - \bxd^{\tpose} \breve{\br}_{\text{im,c}}^{k} + \btx_{\text{du}}^{\tpose} \br_{\text{im,c}}^{k},\\
		&= - \bxd^{\tpose} \breve{\br}_{\text{im,c}}^{k} + \btx_{\text{du}}^{\tpose} \br_{\text{im,c}}^{k} +  \btxdu^{\tpose} \breve{\br}_{\text{im,c}}^{k} - \btxdu^{\tpose} \breve{\br}_{\text{im,c}}^{k},\\
		&= -\big(\bxd - \btxdu\big)^{\tpose} \breve{\br}_{\text{im,c}}^{k} + \btxdu^{\tpose} \big(\br_{\text{im,c}}^{k} - \breve{\br}_{\text{im,c}}^{k}\big).
		\end{aligned}
		\end{equation}
		Subsequent to this, we use the expression for the dual system \cref{eq:fom_im1_dual}  and its residual \cref{eq:resddual_im1_corr} we obtain
		\begin{equation}
		\begin{aligned}
		\brdu &= \bE_{\text{du}} \bxd - \bE_{\text{du}} \btxdu = \bE_{\text{du}}\big( \bxd - \btxdu\big),\\
		\implies & \big( \bxd - \btxdu\big) = \bE_{\text{du}}^{-1} \brdu = \bE_{\text{im}}^{-\tpose} \brdu.
		\end{aligned}
		\label{eq:err_resd_reln}
		\end{equation}
		Substituting \cref{eq:err_resd_reln} into \cref{eq:t_ymod2} yields
		\begin{align}
		\by^{k} - \overline{\by}_{\text{im,c}}^{k} = - \brdu^{\tpose} \bE_{\text{im}}^{-1} \breve{\br}_{\text{im,c}}^{k} + \btxdu^{\tpose} \big(\br_{\text{im,c}}^{k} - \breve{\br}_{\text{im,c}}^{k}\big).
		\end{align}
		Taking the norm on either sides and using the triangle and Cauchy-Schwartz inequalities we obtain the error bound as
		\begin{align}
		\nrm{\by^{k} - \overline{\by}_{\text{im,c}}^{k}} &= \|- \brdu^{\tpose} \bE_{\text{im}}^{-1} \breve{\br}_{\text{im,c}}^{k} \| + \| \btxdu^{\tpose} \big(\br_{\text{im,c}}^{k} - \breve{\br}_{\text{im,c}}^{k}\big) \|,\\
		&\leq \| \bE_{\text{im}}^{-1} \| \| \brdu \| \| \breve{\br}_{\text{im,c}}^{k}\| + \| \btxdu \| \| \br_{\text{im,c}}^{k} - \breve{\br}_{\text{im,c}}^{k}\big) \|.
		\end{align}
	\end{proof}%
\end{appendices}

\addcontentsline{toc}{section}{References}
\bibliographystyle{plainurl}
\bibliography{refs}

\begin{thebibliography}{10}

\bibitem{AbhyankarEtAl2018}
Shrirang Abhyankar, Jed Brown, Emil~M. Constantinescu, Debojyoti Ghosh,
  Barry~F. Smith, and Hong Zhang.
\newblock {PETSc/TS}: {A} modern scalable {ODE/DAE} solver library.
\newblock Technical report, 2018.
\newblock \href {http://arxiv.org/abs/1806.01437} {\path{arXiv:1806.01437}}.

\bibitem{AscRW95}
Uri~M. Ascher, Steven~J. Ruuth, and Brian T.~R. Wetton.
\newblock Implicit-explicit methods for time-dependent partial differential
  equations.
\newblock {\em SIAM J. Numer. Anal.}, 32(3):797--823, 1995.
\newblock \href {https://doi.org/10.1137/0732037} {\path{doi:10.1137/0732037}}.

\bibitem{petsc-web-page}
Satish Balay, Shrirang Abhyankar, Mark~F. Adams, Steven Benson, Jed Brown,
  Peter Brune, Kris Buschelman, Emil~M. Constantinescu, Lisandro Dalcin, Alp
  Dener, Victor Eijkhout, Jacob Faibussowitsch, William~D. Gropp, V\'{a}clav
  Hapla, Tobin Isaac, Pierre Jolivet, Dmitry Karpeev, Dinesh Kaushik,
  Matthew~G. Knepley, Fande Kong, Scott Kruger, Dave~A. May, Lois~Curfman
  McInnes, Richard~Tran Mills, Lawrence Mitchell, Todd Munson, Jose~E. Roman,
  Karl Rupp, Patrick Sanan, Jason Sarich, Barry~F. Smith, Stefano Zampini, Hong
  Zhang, Hong Zhang, and Junchao Zhang.
\newblock {PETS}c {W}eb page.
\newblock \url{https://petsc.org/}, 2023.
\newblock URL: \url{https://petsc.org/}.

\bibitem{morBarMNetal04}
M.~Barrault, Y.~Maday, N.~C. Nguyen, and A.~T. Patera.
\newblock An `empirical interpolation' method: application to efficient
  reduced-basis discretization of partial differential equations.
\newblock {\em C.R. Acad. Sci. Paris}, 339(9):667--672, 2004.
\newblock \href {https://doi.org/10.1016/j.crma.2004.08.006}
  {\path{doi:10.1016/j.crma.2004.08.006}}.

\bibitem{morBenetala21}
P.~Benner, S.~Grivet-Talocia, A.~Quarteroni, G.~Rozza, and L.~M. Schilder,
  W.~Silveira, editors.
\newblock {\em {M}odel {O}rder {R}eduction. Volume 1: {S}ystem- and
  {D}ata-{D}riven {M}ethods and {A}lgorithms}.
\newblock De Gruyter, 2021.
\newblock \href {https://doi.org/10.1515/9783110499001}
  {\path{doi:10.1515/9783110499001}}.

\bibitem{morBenetalb21}
P.~Benner, S.~Grivet-Talocia, A.~Quarteroni, G.~Rozza, and L.~M. Schilder,
  W.~Silveira, editors.
\newblock {\em {M}odel {O}rder {R}eduction. Volume 2: {S}napshot-{B}ased
  {M}ethods and {A}lgorithms}.
\newblock De Gruyter, 2021.
\newblock \href {https://doi.org/10.1515/9783110671490}
  {\path{doi:10.1515/9783110671490}}.

\bibitem{morBenetalc21}
P.~Benner, S.~Grivet-Talocia, A.~Quarteroni, G.~Rozza, and L.~M. Schilder,
  W.~Silveira, editors.
\newblock {\em {M}odel {O}rder {R}eduction. Volume 3: Applications}.
\newblock De Gruyter, 2021.
\newblock \href {https://doi.org/10.1515/9783110499001}
  {\path{doi:10.1515/9783110499001}}.

\bibitem{Bis06}
Christopher~M. Bishop.
\newblock {\em Pattern recognition and machine learning}.
\newblock Information Science and Statistics. Springer, New York, 2006.
\newblock \href {https://doi.org/10.1007/978-0-387-45528-0}
  {\path{doi:10.1007/978-0-387-45528-0}}.

\bibitem{morBufetal12}
Annalisa Buffa, Yvon Maday, Anthony~T. Patera, Christophe Prud'homme, and
  Gabriel Turinici.
\newblock {{A} priori} convergence of the greedy algorithm for the parametrized
  reduced basis method.
\newblock {\em ESAIM Math. Model. Numer. Anal.}, 46(3):595--603, 2012.
\newblock \href {https://doi.org/10.1051/m2an/2011056}
  {\path{doi:10.1051/m2an/2011056}}.

\bibitem{Buh03}
M.~D. Buhmann.
\newblock {\em Radial {B}asis {F}unctions: {T}heory and {I}mplementations},
  volume~12 of {\em Cambridge Monographs on Applied and Computational
  Mathematics}.
\newblock Cambridge University Press, Cambridge, 2003.

\bibitem{But08}
J.~C. Butcher.
\newblock {\em Numerical methods for ordinary differential equations}.
\newblock John Wiley \& Sons, Ltd., Chichester, {S}econd edition, 2008.
\newblock \href {https://doi.org/10.1002/9780470753767}
  {\path{doi:10.1002/9780470753767}}.

\bibitem{morCanTU09}
C.~Canuto, T.~Tonn, and K.~Urban.
\newblock A posteriori error analysis of the reduced basis method for nonaffine
  parametrized nonlinear {PDE}s.
\newblock {\em {SIAM} J. Numer. Anal.}, 47(3):2001--2022, 2009.
\newblock \href {https://doi.org/10.1137/080724812}
  {\path{doi:10.1137/080724812}}.

\bibitem{morCarBF11}
K.~Carlberg, C.~Bou-Mosleh, and C.~Farhat.
\newblock Efficient non-linear model reduction via a least-squares
  {P}etrov-{G}alerkin projection and compressive tensor approximations.
\newblock {\em Internat. J. Numer. Methods Engrg.}, 86(2):155--181, 2011.
\newblock \href {https://doi.org/10.1002/nme.3050}
  {\path{doi:10.1002/nme.3050}}.

\bibitem{morCasetal15}
Fabien Casenave, Alexandre Ern, and Tony Leli\`evre.
\newblock A nonintrusive reduced basis method applied to aeroacoustic
  simulations.
\newblock {\em Adv. Comput. Math.}, 41(5):961--986, 2015.
\newblock \href {https://doi.org/10.1007/s10444-014-9365-0}
  {\path{doi:10.1007/s10444-014-9365-0}}.

\bibitem{morChaH18}
R.~Chakir and J.K. Hammond.
\newblock A non-intrusive reduced basis method for elastoplasticity problems in
  geotechnics.
\newblock {\em J. Comput. Appl. Math.}, 337:1--17, 2018.
\newblock \href {https://doi.org/10.1016/j.cam.2017.12.044}
  {\path{doi:10.1016/j.cam.2017.12.044}}.

\bibitem{morChaM09}
Rachida Chakir and Yvon Maday.
\newblock {A two-grid finite-element/reduced basis scheme for the approximation
  of the solution of parameter dependent PDE}.
\newblock In {\em {9e Colloque National en Calcul des Structures}}, Giens,
  France, 2009. {CSMA}.
\newblock URL: \url{https://hal.archives-ouvertes.fr/hal-01420726}.

\bibitem{morChaS10}
S.~Chaturantabut and D.~C. Sorensen.
\newblock Nonlinear model reduction via discrete empirical interpolation.
\newblock {\em {SIAM} J. Sci. Comput.}, 32(5):2737--2764, 2010.
\newblock \href {https://doi.org/10.1137/090766498}
  {\path{doi:10.1137/090766498}}.

\bibitem{morCheFB19a}
S.~Chellappa, L.~Feng, and P.~Benner.
\newblock Adaptive basis construction and improved error estimation for
  parametric nonlinear dynamical systems.
\newblock {\em Internat. J. Numer. Methods Engrg.}, 121(23):5320--5349, 2020.
\newblock \href {https://doi.org/10.1002/nme.6462}
  {\path{doi:10.1002/nme.6462}}.

\bibitem{morChe23}
Sridhar Chellappa.
\newblock {\em A Posteriori Error Estimation and Adaptivity for Model Order
  Reduction of Large-Scale Systems}.
\newblock {D}issertation, Otto-von-Guericke-Universit{\"a}t, Magdeburg,
  Germany, 2023.
\newblock \href {https://doi.org/http://dx.doi.org/10.25673/101396}
  {\path{doi:http://dx.doi.org/10.25673/101396}}.

\bibitem{morCheQR17}
Peng Chen, Alfio Quarteroni, and Gianluigi Rozza.
\newblock Reduced basis methods for uncertainty quantification.
\newblock {\em SIAM/ASA J. Uncertain. Quantif.}, 5(1):813--869, 2017.
\newblock \href {https://doi.org/10.1137/151004550}
  {\path{doi:10.1137/151004550}}.

\bibitem{Cheetal18}
Ricky T.~Q. Chen, Yulia Rubanova, Jesse Bettencourt, and David~K Duvenaud.
\newblock Neural ordinary differential equations.
\newblock In S.~Bengio, H.~Wallach, H.~Larochelle, K.~Grauman, N.~Cesa-Bianchi,
  and R.~Garnett, editors, {\em Advances in Neural Information Processing
  Systems}, volume~31. Curran Associates, Inc., 2018.
\newblock URL:
  \url{https://proceedings.neurips.cc/paper/2018/file/69386f6bb1dfed68692a24c8686939b9-Paper.pdf}.

\bibitem{CohD15}
Albert Cohen and Ronald DeVore.
\newblock {Kolmogorov widths under holomorphic mappings}.
\newblock {\em IMA Journal of Numerical Analysis}, 36(1):1--12, 03 2015.
\newblock \href {https://doi.org/10.1093/imanum/dru066}
  {\path{doi:10.1093/imanum/dru066}}.

\bibitem{morDegVW20}
Denise Degen, Karen Veroy, and Florian Wellmann.
\newblock Certified reduced basis method in geosciences: addressing the
  challenge of high-dimensional problems.
\newblock {\em Comput. Geosci.}, 24(1):241--259, 2020.
\newblock \href {https://doi.org/10.1007/s10596-019-09916-6}
  {\path{doi:10.1007/s10596-019-09916-6}}.

\bibitem{morDihH16}
Markus Dihlmann and Bernard Haasdonk.
\newblock A reduced basis {K}alman filter for parametrized partial differential
  equations.
\newblock {\em ESAIM Control Optim. Calc. Var.}, 22(3):625--669, 2016.
\newblock \href {https://doi.org/10.1051/cocv/2015019}
  {\path{doi:10.1051/cocv/2015019}}.

\bibitem{morDroHO12}
M.~Drohmann, B.~Haasdonk, and M.~Ohlberger.
\newblock Reduced basis approximation for nonlinear parametrized evolution
  equations based on empirical operator interpolation.
\newblock {\em {SIAM} J. Sci. Comput.}, 34(2):A937--A969, 2012.
\newblock \href {https://doi.org/10.1137/10081157X}
  {\path{doi:10.1137/10081157X}}.

\bibitem{gardner2022sundials}
David~J Gardner, Daniel~R Reynolds, Carol~S Woodward, and Cody~J Balos.
\newblock Enabling new flexibility in the {SUNDIALS} suite of nonlinear and
  differential/algebraic equation solvers.
\newblock {\em ACM Transactions on Mathematical Software (TOMS)}, 2022.
\newblock \href {https://doi.org/10.1145/3539801} {\path{doi:10.1145/3539801}}.

\bibitem{GreU19}
Constantin Greif and Karsten Urban.
\newblock Decay of the {K}olmogorov {$N$}-width for wave problems.
\newblock {\em Appl. Math. Lett.}, 96:216--222, 2019.
\newblock \href {https://doi.org/10.1016/j.aml.2019.05.013}
  {\path{doi:10.1016/j.aml.2019.05.013}}.

\bibitem{morGre05}
M.~Grepl.
\newblock {\em Reduced-basis approximation a posteriori error estimation for
  parabolic partial differential equations}.
\newblock PhD thesis, Massachussetts Institute of Technology (MIT), Cambridge,
  USA, 2005.
\newblock URL: \url{http://dspace.mit.edu/handle/1721.1/7582}.

\bibitem{morGreP05}
M.~A. Grepl and A.~T. Patera.
\newblock A posteriori error bounds for reduced-basis approximations of
  parametrized parabolic partial differential equations.
\newblock {\em {ESAIM}: Math. Model. Numer. Anal.}, 39(1):157--181, 2005.
\newblock \href {https://doi.org/10.1051/m2an:2005006}
  {\path{doi:10.1051/m2an:2005006}}.

\bibitem{morGre12}
Martin~A. Grepl.
\newblock Certified reduced basis methods for nonaffine linear time-varying and
  nonlinear parabolic partial differential equations.
\newblock {\em Math. Models Methods Appl. Sci.}, 22(3):1150015, 40, 2012.
\newblock \href {https://doi.org/10.1142/S0218202511500151}
  {\path{doi:10.1142/S0218202511500151}}.

\bibitem{morGroM22}
Elise Grosjean and Yvon Maday.
\newblock Error estimate of the non-intrusive reduced basis {(NIRB)} two-grid
  method with parabolic equations.
\newblock e-prints 2211.08897, arXiv, 2022.
\newblock URL: \url{https://arxiv.org/abs/2211.08897}, \href
  {https://doi.org/10.48550/arXiv.2211.08897}
  {\path{doi:10.48550/arXiv.2211.08897}}.

\bibitem{morHaaO08}
B.~Haasdonk and M.~Ohlberger.
\newblock Reduced basis method for finite volume approximations of parametrized
  linear evolution equations.
\newblock {\em {ESAIM}: Math. Model. Numer. Anal.}, 42(2):277 -- 302, 2008.
\newblock \href {https://doi.org/10.1051/m2an:2008001}
  {\path{doi:10.1051/m2an:2008001}}.

\bibitem{morHaaO11}
B.~Haasdonk and M.~Ohlberger.
\newblock Efficient reduced models and a posteriori error estimation for
  parametrized dynamical systems by offline/online decomposition.
\newblock {\em Math. Comput. Model. Dyn. Syst.}, 17(2):145--161, 2011.
\newblock \href {https://doi.org/10.1080/13873954.2010.514703}
  {\path{doi:10.1080/13873954.2010.514703}}.

\bibitem{morharthw18}
Dirk Hartmann, Matthias Herz, and Utz Wever.
\newblock {\em Model Order Reduction a Key Technology for Digital Twins}, pages
  167--179.
\newblock Springer-Verlag, Cham, 2018.
\newblock \href {https://doi.org/10.1007/978-3-319-75319-5_8}
  {\path{doi:10.1007/978-3-319-75319-5_8}}.

\bibitem{morHesRS16}
J.~S. Hesthaven, G.~Rozza, and B.~Stamm.
\newblock {\em Certified Reduced Basis Methods for Parametrized Partial
  Differential Equations}.
\newblock SpringerBriefs in Mathematics. Springer International Publishing,
  2016.
\newblock \href {https://doi.org/10.1007/978-3-319-22470-1}
  {\path{doi:10.1007/978-3-319-22470-1}}.

\bibitem{hindmarsh1983odepack}
Alan~C Hindmarsh.
\newblock {ODEPACK}, a systematized collection of ode solvers.
\newblock In R.~S. Stepleman, editor, {\em Scientific computing : applications
  of mathematics and computing to the physical sciences}. Elsevier, 1983.

\bibitem{hindmarsh2005sundials}
Alan~C Hindmarsh, Peter~N Brown, Keith~E Grant, Steven~L Lee, Radu Serban,
  Dan~E Shumaker, and Carol~S Woodward.
\newblock {SUNDIALS}: Suite of nonlinear and differential/algebraic equation
  solvers.
\newblock {\em ACM Transactions on Mathematical Software (TOMS)},
  31(3):363--396, 2005.
\newblock \href {https://doi.org/10.1145/1089014.1089020}
  {\path{doi:10.1145/1089014.1089020}}.

\bibitem{hinesrbf}
Trever Hines.
\newblock Python package containing tools for radial basis function ({RBF})
  applications.
\newblock \url{https://github.com/treverhines/RBF}, 2023.

\bibitem{Huaetal22}
Zhongzhan Huang, Senwei Liang, Hong Zhang, Haizhao Yang, and Liang Lin.
\newblock Accelerating numerical solvers for large-scale simulation of
  dynamical system via neurvec.
\newblock e-prints 2208.03680, arXiv, 2022.
\newblock cs.CE.
\newblock URL: \url{https://arxiv.org/abs/2208.03680}.

\bibitem{morKapetal22}
M.~G. Kapteyn, D.~J. Knezevic, D.~B.~P. Huynh, M.~Tran, and K.~E. Willcox.
\newblock Data-driven physics-based digital twins via a library of
  component-based reduced-order models.
\newblock {\em Internat. J. Numer. Methods Engrg.}, 123(13):2986--3003, 2022.
\newblock \href {https://doi.org/10.1002/nme.6423}
  {\path{doi:10.1002/nme.6423}}.

\bibitem{morKaretal18}
Mark K\"{a}rcher, S\'{e}bastien Boyaval, Martin~A. Grepl, and Karen Veroy.
\newblock Reduced basis approximation and a posteriori error bounds for
  4{D}-{V}ar data assimilation.
\newblock {\em Optim. Eng.}, 19(3):663--695, 2018.
\newblock \href {https://doi.org/10.1007/s11081-018-9389-2}
  {\path{doi:10.1007/s11081-018-9389-2}}.

\bibitem{LeeS19}
Ju~Weon {Lee} and Andreas {Seidel-Morgenstern}.
\newblock {Solving hyperbolic conservation laws with active counteraction
  against numerical errors: Isothermal fixed-bed adsorption}.
\newblock {\em Chemical Engineering Science}, 207:1309--1330, 2019.
\newblock \href {https://doi.org/10.1016/j.ces.2019.07.053}
  {\path{doi:10.1016/j.ces.2019.07.053}}.

\bibitem{morMacMP01}
L.~Machiels, Y.~Maday, and A.~T. Patera.
\newblock Output bounds for reduced-order approximations of elliptic partial
  differential equations.
\newblock {\em Comp. Meth. Appl. Mech. Eng.}, 190(26-27):3413--3426, 2001.
\newblock \href {https://doi.org/10.1016/S0045-7825(00)00275-9}
  {\path{doi:10.1016/S0045-7825(00)00275-9}}.

\bibitem{Pin85}
Allan Pinkus.
\newblock {\em {$n$}-widths in approximation theory}, volume~7 of {\em
  Ergebnisse der Mathematik und ihrer Grenzgebiete (3) [Results in Mathematics
  and Related Areas (3)]}.
\newblock Springer-Verlag, Berlin, 1985.
\newblock \href {https://doi.org/10.1007/978-3-642-69894-1}
  {\path{doi:10.1007/978-3-642-69894-1}}.

\bibitem{Polietal20}
Michael Poli, Stefano Massaroli, Atsushi Yamashita, Hajime Asama, and Jinkyoo
  Park.
\newblock Hypersolvers: Toward fast continuous-depth models.
\newblock In H.~Larochelle, M.~Ranzato, R.~Hadsell, M.F. Balcan, and H.~Lin,
  editors, {\em Advances in Neural Information Processing Systems}, volume~33,
  pages 21105--21117. Curran Associates, Inc., 2020.
\newblock URL:
  \url{https://proceedings.neurips.cc/paper/2020/file/f1686b4badcf28d33ed632036c7ab0b8-Paper.pdf}.

\bibitem{morQuaMN16}
A.~Quarteroni, A.~Manzoni, and F.~Negri.
\newblock {\em {R}educed {B}asis {M}ethods for {P}artial {D}ifferential
  {E}quations}, volume~92 of {\em La Matematica per il 3+2}.
\newblock Springer International Publishing, 2016.
\newblock \href {https://doi.org/10.1007/978-3-319-15431-2}
  {\path{doi:10.1007/978-3-319-15431-2}}.

\bibitem{reynolds2022arkode}
Daniel~R Reynolds, David~J Gardner, Carol~S Woodward, and Rujeko Chinomona.
\newblock {ARKODE: A flexible IVP solver infrastructure for one-step methods}.
\newblock {\em arXiv preprint arXiv:2205.14077}, 2022.

\bibitem{morRov03}
Dimitrios~V. Rovas.
\newblock {\em Reduced-Basis Output Bound Methods for Parametrized Partial
  Differential Equations}.
\newblock PhD thesis, Massachussetts Institute of Technology (MIT), Cambridge,
  USA, 2003.
\newblock URL: \url{https://dspace.mit.edu/handle/1721.1/16956}.

\bibitem{morRozHP08}
G.~Rozza, D.~B.~P. Huynh, and A.~T. Patera.
\newblock Reduced basis approximation and a posteriori error estimation for
  affinely parametrized elliptic coercive partial differential equations:
  application to transport and continuum mechanics.
\newblock {\em Arch. Comput. Methods Eng.}, 15(3):229--275, 2008.
\newblock \href {https://doi.org/10.1007/s11831-008-9019-9}
  {\path{doi:10.1007/s11831-008-9019-9}}.

\bibitem{matlabodesuite}
Lawrence~F. Shampine and Mark~W. Reichelt.
\newblock The {MATLAB} {ODE} suite.
\newblock volume~18, pages 1--22. 1997.
\newblock Dedicated to C. William Gear on the occasion of his 60th birthday.
\newblock \href {https://doi.org/10.1137/S1064827594276424}
  {\path{doi:10.1137/S1064827594276424}}.

\bibitem{SheCL20}
Xing Shen, Xiaoliang Cheng, and Kewei Liang.
\newblock Deep euler method: solving odes by approximating the local truncation
  error of the euler method.
\newblock e-prints 2003.09573, arXiv, 2020.
\newblock URL: \url{https://arxiv.org/abs/2003.09573}, \href
  {https://doi.org/10.48550/arXiv.2003.09573}
  {\path{doi:10.48550/arXiv.2003.09573}}.

\bibitem{morVeretal03}
K.~Veroy, C.~Prud'Homme, D.~V. Rovas, and A.~T. Patera.
\newblock A posteriori error bounds for reduced-basis approximation of
  parametrized noncoercive and nonlinear elliptic partial differential
  equations.
\newblock In {\em {16th AIAA Computational Fluid Dynamics Conference}},
  Orlando, United States, 2003.
\newblock URL: \url{https://hal.archives-ouvertes.fr/hal-01219051}.

\bibitem{morWanHR19}
Qian Wang, Jan~S. Hesthaven, and Deep Ray.
\newblock Non-intrusive reduced order modeling of unsteady flows using
  artificial neural networks with application to a combustion problem.
\newblock {\em J. Comput. Phys.}, 384:289--307, 2019.
\newblock \href {https://doi.org/10.1016/j.jcp.2019.01.031}
  {\path{doi:10.1016/j.jcp.2019.01.031}}.

\bibitem{Wed05}
Holger Wendland.
\newblock {\em Scattered {D}ata {A}pproximation}, volume~17 of {\em Cambridge
  Monographs on Applied and Computational Mathematics}.
\newblock Cambridge University Press, Cambridge, 2005.

\bibitem{morWirSH14}
D.~Wirtz, D.~C. Sorensen, and B.~Haasdonk.
\newblock A posteriori error estimation for {DEIM} reduced nonlinear dynamical
  systems.
\newblock {\em {SIAM} J. Sci. Comput.}, 36(2):A311--A338, 2014.
\newblock \href {https://doi.org/10.1137/120899042}
  {\path{doi:10.1137/120899042}}.

\bibitem{morZhaFLetal14}
Y.~Zhang, L.~Feng, S.~Li, and P.~Benner.
\newblock Accelerating {PDE} constrained optimization by the reduced basis
  method: application to batch chromatography.
\newblock {\em Internat. J. Numer. Methods Engrg.}, 104(11):983--1007, 2015.
\newblock \href {https://doi.org/10.1002/nme.4950}
  {\path{doi:10.1002/nme.4950}}.

\bibitem{morZhaFLetal15}
Y.~Zhang, L.~Feng, S.~Li, and P.~Benner.
\newblock An efficient output error estimation for model order reduction of
  parametrized evolution equations.
\newblock {\em {SIAM} J. Sci. Comput.}, 37(6):B910--B936, 2015.
\newblock \href {https://doi.org/10.1137/140998603}
  {\path{doi:10.1137/140998603}}.

\end{thebibliography}
\end{document}